\documentclass[12pt,a4paper]{article}
\usepackage{amsfonts}
\usepackage{amssymb}
\usepackage{amsmath}
\usepackage{latexsym}
\usepackage{graphicx}
\usepackage{pbez}

\def\call#1{\ensuremath{\cal #1}}
\def\M#1{\ensuremath{\mathbb #1}}
\def\f{\ensuremath{\varphi}}

\def\wt#1{\ensuremath{\widetilde{#1}}}
\def\z{{\bf z}}
\def\wb#1{\ensuremath{\overline{#1}}}
\newtheorem{teo}{Theorem}[section]
\newtheorem{prop}[teo]{Proposition}

\newtheorem{lemma}[teo]{Lemma}
\newtheorem{defi}[teo]{Definition}
\newtheorem{remark}[teo]{Remark}

\newcommand{\proof}{\noindent{\em Proof.}\ }
\newcommand{\qed}{\\\makebox[\textwidth][r]{$\Box$}\vspace{\baselineskip}}

\begin{document}
\begin{center}
  {\large Algebraic and geometric solutions of hyperbolic Dehn filling
  equations.}

\vskip5mm
\textsc{Stefano Francaviglia\footnote{Author address: Stefano
      Francaviglia, Scuola Normale Superiore, Piazza dei Cavalieri 7,
      I-56126 Pisa, e-mail: s.francaviglia@sns.it.}}

\today
\end{center}

\vskip5mm

\begin{abstract}
In this paper we study the difference between
algebraic and geometric solutions of the hyperbolic Dehn filling
equations (see Definition~\ref{d as}) for ideally triangulated 3-manifolds.

We show that any geometric solution is an algebraic one, 
and we prove the uniqueness of the geometric solutions.

Then we do explicit calculations for three interesting
examples. With the first two examples we see 
that not all algebraic solutions are geometric and that the algebraic
solutions are not unique. 

The third example is a non-hyperbolic
manifold that admits a positive, partially flat solution of the
compatibility and completeness equations.
\end{abstract}

\section{Introduction}

One of the most useful tools for studying the hyperbolic structures on
3-manifolds is the technique of ideal triangulations,
introduced by Thurston in \cite{Thu} to study the hyperbolic structure
of the complement of the figure-eight knot.
An ideal triangulation of an open 3-manifold $M$ is a description of
$M$ as a disjoint union of copies of the standard tetrahedron with
vertices removed (ideal tetrahedron), glued together by a given set of
face-pairing maps.

Given an ideal triangulation on an
open 3-manifold $M$ with toroidal ends (this is known to be necessary
for hyperbolicity)
the idea is to construct a hyperbolic structure on $M$ by defining it
on each tetrahedron and then by requiring that such structures are
compatible with a global one on $M$. See \cite{Thu}, \cite{NZ},
\cite{BP} for more details.

A complete finite-volume hyperbolic structure with totally
geodesic faces on an oriented ideal tetrahedron is described by
a complex number with positive imaginary part, called modulus (see
Section~\ref{s itwm}).

The compatibility conditions on the hyperbolic structures of
the tetrahedra translate to algebraic equations on the moduli.
These equations depend on the combinatorics of the triangulation and
are called {\em compatibility equations} (see Section~\ref{s itwm}).
When the moduli induce a structure on $M$,
one can ask about the completeness of such a structure, and this
translates to other algebraic equations on the moduli, called
{\em completeness equations}.
Once a complete structure on $M$ is found, one can study small
perturbations by taking moduli which are near to the
complete solution and which satisfy the compatibility equations.
By studying the completions of the structures obtained in this way,
one can see that in certain cases one obtains a
Dehn filling of $M$ for suitable filling-parameters. Conversely, given a
Dehn filling $N$ of $M$ with fixed
filling-parameters, then there exists a system of equations on the moduli,
which express the fact that
the completion of the hyperbolic structure
induced on $M$ by the moduli is exactly $N$.  We call such equations
{\em hyperbolic Dehn filling equations}. The completeness equations
can be viewed as a particular case of hyperbolic Dehn filling
equations, relative to the empty filling.

As the geometric conditions translate to algebraic equations, the
problem can be studied from an algebraic point of view. 
The question that naturally  arises is

\vskip1.5ex
\noindent{\bf Question 1} {\em Does any solution of the algebraic equations
    have a geometric meaning?}

\vskip1.5ex
It is well-known that when all the moduli have a positive imaginary part,
the solutions of the compatibility equations have a
natural geometric interpretation as a decomposition of $M$ into
geodesic ideal tetrahedra.
In general it is not clear when an algebraic solution of the
equations have a geometric interpretation leading to define a
finite-volume hyperbolic structure on $M$.
Moreover, it is still unknown whether a hyperbolic manifolds admits
a decomposition into ``positive'' geodesic ideal tetrahedra.

Epstein and Penner
\cite{ep} have shown that a decomposition into convex geodesic ideal
polyhedra always exists. By subdividing such a  decomposition flat
tetrahedra can appear. This translates to the fact that some modulus
becomes real. Petronio and Weeks \cite{PW} have shown that a partially
flat solution of the compatibility and completeness equations leads to
a complete hyperbolic structure on $M$, while a solution of the
compatibility equations alone does not in general lead to an (even
incomplete) hyperbolic structure 
(we notice that in~\cite{PW} the system of compatibility
equations contains additional equations on the angles of the moduli). 
By perturbing a partially flat
decomposition of $M$, negatively oriented tetrahedra appear.
This translates to the fact that the imaginary part of some modulus
becomes negative. The geometric meaning of a solution of the equations,
that involves positive, negative and flat tetrahedra is not clear.
Petronio and Porti \cite{PP} have shown that 
the results known in the case in which all the tetrahedra are positive 
extend near a partially flat triangulation obtained as a subdivision of an
Epstein-Penner decomposition. This leads to the following problem.

\vskip1.5ex
\noindent{\bf Question 2} {\em To understand the geometric meaning of
  the mixed solutions.}

\vskip1.5ex
Finally one can ask about uniqueness

\vskip1.5ex
\noindent{\bf Question 3} {\em Are the algebraic/geometric solutions
  unique?}

\vskip2.5ex
\noindent
In this paper we concentrate on the differences between
algebraic and geometric solutions of hyperbolic Dehn filling
equations. 
Our interpretation of mixed solutions will be in terms of developing
maps and holonomy. Roughly speaking, we will call {\em geometric solution} a
choice of moduli whose holonomy is well-defined and discrete 
(see Definition~\ref{d as} for the exact definition).

In this setting, we show that {\em algebraic is strictly bigger than
  geometric},in the sense that any geometric solution is algebraic but
  non all algebraic solutions are geometric, answering to Questions 1
  and 2. We remark that the difference between algebraic and
  geometric solutions appears only with the presence of negative
  tetrahedra. 

Then we show that the geometric solutions are unique while the
algebraic ones are not, answering to Question 3.

The paper is structured as follows:

In Section~\ref{s itwm} we define the notions of developing map and
holonomy for a triangulation with moduli.

In Section~\ref{s hde} we describe the system of hyperbolic Dehn
filling equations, we get the definition of geometric solution and we
prove the following fact:

\vskip1.5ex
\noindent{\bf Theorem~\ref{t gisa}} {\em  Each geometric solution of
  the hyperbolic Dehn filling equations is also an algebraic one.}

\vskip1.5ex
In Section~\ref{s 4} we prove the uniqueness of the geometric solutions.

\vskip1.5ex
\noindent{\bf Theorem~\ref{p uz}} {\em There exists at most one geometric
  solution of the hyperbolic Dehn filling equations.}

\vskip1.5ex
In Section~\ref{s 5} we do explicit calculations
for some interesting examples:
\begin{itemize}
  \item First we study two one-cusped manifolds, that are bundles over
  $S^1$ with fiber a punctured torus and are called $LR^3$ and
  $L^2R^3$. These manifolds admit more algebraic solutions and the
  unique geometric one.
\item Then we study a manifold with non-trivial JSJ decomposition,
  obtained by gluing a Seifert manifold to the complement of the
  figure-eight knot. 
  This manifold is not hyperbolic but it admits a partially flat solution
  of compatibility and completeness equations. 
\end{itemize}

The manifold $LR^3$ (and $L^2R^3$)
 is interesting because, on one hand it shows that the
algebraic solutions are not unique, on the other hand it provides an
example of an algebraic solution which is not geometric
(Proposition~\ref{p nd}).
We notice that such a ``bad'' solution does not involve flat
tetrahedra and has a good behavior on the boundary. Namely, the
boundary torus inherits a intrinsic Euclidean structure (up to scaling).
This fact is surprising because the geometry of a finite-volume
hyperbolic 3-manifold is strictly related to the geometry of its
boundary.

Actually the equations on the moduli have an
interpretation as conditions on the geometry of the
boundary. More precisely, any
ideal triangulation of $M$ induces a triangulation of the boundary
tori, by considering the manifold with boundary obtained by chopping
off an open regular neighborhood of the ideal vertices.
A modulus for the hyperbolic structure of an ideal tetrahedron
determines a modulus for the similarity structure of the triangles
obtained by horospherical sections near the vertices. So an ideal
triangulation with moduli of $M$ induces a triangulation with moduli
of the boundary tori. The compatibility equations
express the fact that the moduli for the triangles lead to
similarity structures on the tori. The completeness equations express
the fact that the structures of the boundary tori are Euclidean.

Moreover, when the imaginary part of the moduli is not negative, the
control of what happens on the boundary it suffices to guarantee the
structure on the whole $M$. For example, in order to
have a complete finite-volume hyperbolic structure on $M$, it suffices
to check that the boundary tori have Euclidean structures.

In \cite{F} it is shown that any algebraic solution of the
compatibility and completeness
equations for the similarity structure of a triangulated torus leads to
a Euclidean structure, even if there are negative triangles,
provided that the algebraic sum of the areas of the triangles is not zero.
So the example of $LR^3$ shows that the Euclidean situation in dimension 2 and
the hyperbolic one in dimension 3 become quite different when we allow
the moduli to have negative imaginary part.

Finally, the manifold with non-trivial JSJ decomposition that we study in 
the last example 
is a manifold that admits a positive, partially flat triangulation of the
compatibility and completeness equations. Such a solution cannot be
geometric as the manifold is not hyperbolic. This seems to
contradicts~\cite{PW}. Actually there is no contradictions because in
our example the conditions on the angles are not satisfied. We notice
that in general  not any geometric solution satisfies the conditions
on the angles. Nevertheless, this example shows that such conditions
play a central role for a solution to be geometric.

\section{Ideal triangulations with moduli, developing maps and
  holonomy}\label{s itwm}

Let $M$ be the interior of a compact 3-manifold $\wb M$  with boundary
and let $\widehat M$ be the space obtained from
$\wb M$ by collapsing each component of its boundary to a point. The
space $\widehat M$
is homeomorphic to the space obtained from $\wb M$ by gluing to each
boundary component $C$ the cone over $C$. In the sequel we will often
identify $M$ with its image under the projection $\wb M\to \widehat M$.

Let $\Delta$ be the standard 3-simplex and let $\Delta^*$ be the standard
ideal 3-simplex, i.e. $\Delta $ with vertices removed.

We define now what we mean by ideal triangulation. Roughly speaking,
an ideal triangulation of $M$ is a presentation of $M$ as a union of
ideal tetrahedra. By a face-pairing rule  $r:F_1 \dashrightarrow F_2$
between two-dimensional faces of two tetrahedra (possibly the same
tetrahedron) we mean a bijective correspondence $r$ from 
the vertices of $F_1$ to the ones of $F_2$. A realization of a 
face-pairing rule is a homeomorphism that extends the rule.

\begin{defi}\label{d it}
 Let $\{\Delta_i,\ i\in I\}$ be a set of copies of $\Delta$ and
let $\{r_j:F_{j_1}\dashrightarrow F_{j_2},\ j\in J\}$ be a set of face-pairing
rules between two-dimensional faces of the $\Delta_i$'s. We say that
$\tau=(\{\Delta_i\},\{r_j\})$ is an {\em ideal triangulation} of
  $M$ if for each $j\in J$
there exists a set $\{f_j:F_{j_1}\to F_{j_2},\  j\in J\}$ of
simplicial maps such that each 
$f_j$ is an extension of  the rule $r_j$ and such that there exists
a homeomorphism $\f:(\sqcup_i \Delta_i^*\big)/\{f_j\}\to M$.
If $M$ is oriented, we fix an orientation for
$\Delta$ and we require the $r_j$'s to be orientation reversing and
\f\  to be orientation preserving.
We say that $\call{R}$ $=(\{f_j\},\f)$ is a {\em realization} of $\tau$.

Given a realization of $\tau$, for each $i\in I$ we set
$\f_{_{\Delta_i}}$ to be the composition of \f\ with the inclusion
$\Delta_i^*\to \sqcup_i (\Delta_i^*\big)/\{f_j\}$.
\end{defi}

To require $\f$ to be a homeomorphism from $(\sqcup_i
\Delta_i^*)/\{f_j\}$ to $M$
is equivalent to require $\f$ to extend to a homeomorphism from $(\sqcup_i
\Delta_i)/\{f_j\}$ to $\widehat M$ such that it is one-to-one
from  $(\sqcup_i \Delta_i/\{f_j\})\setminus(\sqcup_i \Delta_i^*/\{f_j\})$ to
$\widehat M \setminus M$. Hence the vertices of $\sqcup_i
\Delta_i/\{f_j\}$ are in one-to-one correspondence with the
connected components of $\partial \wb M$.

Let $U$ be an open regular neighborhood of the
0-skeleton of a tetrahedron $\Delta$,
we call $\Delta^-$ the truncated tetrahedron
$\Delta\setminus U$ and we call $\partial^-\Delta=\partial U$.

Given an ideal triangulation $\tau$ it is possible to truncate each
$\Delta_i$  in such a way that
$(\sqcup_i \Delta_i^-/\{f_j\},\sqcup_i \partial^-\Delta_i/\{f_j\})$ is
homeomorphic to $(\wb M,\partial \wb M)$. In other words, each ideal
triangulation induces a triangulation of $\partial \wb M$ (here
triangulation means possibly with auto-adjacencies).

\ \\
We introduce now the notion of modulus of a hyperbolic ideal
tetrahedron. Let $A$ be a straight ideal tetrahedron in $\M H^3$,
i.e. a geodesic 
tetrahedron
in $\overline{\M H^3}$ such that $\partial \M H^3\cap A$ is the
0-skeleton of $A$. Since such a tetrahedron is the convex hull of its
vertices, it is completely determined by them.

The hyperbolic ideal tetrahedra that we consider in the sequel
can be flat, but not degenerate. This means that a tetrahedron can be
completely contained in a 2-plane, but we always require that its
vertices are four distinct points.

An orientation of an abstract tetrahedron is an ordering of its
vertices up to even permutations.
When $A$ is fat (i.e. it is not contained in a 2-plane), the
orientations of $A$ as a subset of $\M H^3$ are in correspondence
with its orientations as an abstract tetrahedron.

Consider now the half-space model $\M C\times \M R^+$ of $\M H^3$ and let
$(v_1,v_2,v_3,v_4)$ be an ordering of the vertices of $A$. Then
the Isom$^+(\M H^3)$-class of $A$ is determined by a complex
number $z$ by mapping $(v_1,v_2,v_3,v_4)$ to $(0,1,\infty,z)$ via an
isometry. Since Isom$^+(\M H^3)$ acts transitively on the set of
triples of points at infinity,
it is always possible to map $(v_1,v_2,v_3)$ to
$(0,1,\infty)$ via an element of Isom$^+(\M H^3)$. Moreover
such an isometry is unique and the
number $z$ is the cross-ratio $[v_1:v_2:v_3:v_4]$.

Since the vertices of $A$ are four distinct points in $\M
C\cup \{\infty\}$, it follows that $z\notin\{0,1,\infty\}$. 
By slicing
$A$ with a horosphere $\M C\times\{t\}$, for a sufficiently large $t$
we obtain a 
triangle with a Euclidean structure. Up to scaling, this structure is
exactly the one of the triangle in \M C with vertices in
$\{0,1,z\}$. It follows that the hyperbolic structure of $A$ is
determined by the similarity structure (Euclidean up to scaling)
of a horospherical triangle at a vertex of $A$.

It is easily checked that as the ordering of the vertices varies in the
same orientation class, then $z$ varies on the set
$\{z,\frac{1}{1-z},1-\frac{1}{z}\}$. This ambiguity can be avoided by
fixing a preferred edge $e$ of $A$ and arranging the vertices
$(v_1,v_2,v_3,v_4)$ in such a way that $e$ joins $v_1$ and $v_3$ (note
that $[v_1:v_2:v_3:v_4]=[v_3:v_4:v_1:v_2]$).  
The number $z$ is called {\em modulus}
of $A$ associated to $e$.  It is easy to see that
at opposite edges is associated the same modulus.

\begin{remark} In the sequel we tacitly assume
that an orientation and an edge for each tetrahedra have been fixed.
\end{remark}

From now on we consider $M$ to be oriented.

\begin{defi}
Let $\tau$ be an ideal triangulation of $M$.
A choice of moduli $\z=\{z_i,\ i\in I\}$ for $\tau$,
is a choice of a complex number $z_i\neq0,1$ for each tetrahedron
$\Delta_i$ of $\tau$.
We write $(\tau,\z)$ to mean an ideal triangulation $\tau$ with a
choice of moduli $\z$ for $\tau$.
\end{defi}

\begin{defi}
  Let $\Delta$ be the standard tetrahedron. We say that a
  continuous map
  $f:\Delta\to\wb{\M H}^3$ is {\em straight} if 
   \begin{itemize} 
    \item $f$ maps the vertices of $\Delta$ to $\partial \wb{\M H}^3$;
    \item for each sub-simplex $F$ of $\Delta$, the image of $F$ is
     contained in the convex hull of the image of its vertices;
    \item for each sub-simplex $F$ of $\Delta$, the map $f_{|F}$ is a
     homeomorphism whenever $f_{|\partial F}$ is a homeomorphism. 
   \end{itemize} 
\end{defi}

\begin{defi}\label{d cm}

Let $z$ be a choice of a modulus for a tetrahedron $\Delta$
and $\xi:\Delta\to\overline{\M H^3}$ be a continuous map.

We say that $\xi$ is {\em compatible} with $z$ if it is straight
and its image is a geodesic ideal tetrahedron of modulus $z$.

\end{defi}

\begin{defi}\label{d gcr}
Let \call R$=(\{f_j\},\f)$ be a realization of an ideal triangulation
with moduli
$(\tau,\z)$ of $M$. Let $\{g_i:\Delta_i\to \overline{\M H^3}\}$ be a
set of maps, each one compatible with the relative $z_i$. For each
$j\in J$  let $\psi_j$ be the unique orientation-preserving isometry
which realizes the face-pairing rule $r_j$ between the corresponding
faces of the hyperbolic tetrahedra
$g_i(\Delta_i)$.

We say that $\{g_i\}$ is compatible with \call R if for each $j$,
called $\Delta_{i_1}$ and $\Delta_{i_2}$ the tetrahedra glued by the
rule $r_j$, we have that
\begin{equation}\label{eq rel}f_j=g_{i_2}^{-1}\circ\psi_j\circ
g_{i_1}.\end{equation}
\end{defi}

\begin{remark}
  Each time we have some covering and some object $o$ which can be lifted
  in some natural way, as usual we call \wt o a lift of $o$. When it
  is clear, we leave
  to the reader the proofs that what we do is independent from the
  chosen lift. In the following $\wt M$ will be the
  universal covering of $M$ and $\wt\tau$ the lift of $\tau$ on \wt M.
\end{remark}

Given $(\tau,\z)$, the idea is to define a hyperbolic structure on $M$
by taking, for each $i$, an ideal hyperbolic straight tetrahedron of modulus
$z_i$ and then by gluing them realizing the rules $r_j$ via
isometries.

In order to succeed in this construction it is easy to see that a
necessary condition is that for each edge $e$ of $\tau$ the product of
moduli around $e$ is $1$. As mentioned above, if $z_i$ is the modulus of
$\Delta_i$ associated to an edge, then changing edge, the modulus change
in the set $\{z_i,\frac{1}{1-z_i},1-\frac{1}{z_i}\}$. It follows that
the condition of the product of moduli around the edges can be written
as a system \call C of algebraic equations on the moduli, having the form
$$\displaystyle{\pm\prod_i z_i^{\alpha_i}(1-z_i)^{\beta_i}=1}$$ where
$\alpha_i, \beta_i\in\M Z$ depend on the triangulation and on the
chosen preferred edges of the tetrahedra. These equations
are called {\em compatibility equations}.

\begin{lemma}
Let $(\tau,\z)$ be an ideal triangulation with moduli on $M$ such that
equations \call C are satisfied.

Then for each realization \call R$=(\{f_j\},\f)$ of $\tau$ there
exists a set of maps $\{g_i:\Delta_i\to\wb{\M H}^3\}$ compatible with \call R.
\end{lemma}

\proof
We define the $g_i$'s recursively on the $n-$skeleta of $\tau$. On the
0-skeleton we define the maps simply looking at the compatibility
with the moduli. Then a set of maps $\psi_j\in$ Isom$^+(\M H^3)$ as in
Definition~\ref{d gcr} is well-defined.

Let $e$ be an edge of a tetrahedron
$\Delta_{i_0}$ with vertices $e_0$ and $e_1$. We set $g_{i_0}$ on
$e$  to be an homeomorphism onto the geodesic between $g_{i_0}(e_0)$
and $g_{i_0}(e_1)$. Now we define the $g_i$'s on the edges glued to
$e$ by the maps $f_j$ using the formula~\ref{eq rel} of
Definition~\ref{d gcr}. Note that since \call C holds this
 is a good definition. We define the $g_i$'s on the others edges in a
similar way.

Once the $g_i$'s are defined on the 1-skeleton there are no problems
to use again the formula~\ref{eq rel} to define them on the 2-skeleton
and there are no obstructions to extend such maps to the 3-cells.\qed

\begin{remark}
From now on  when we speak about an ideal
triangulation, we suppose that a realization \call R has been fixed
and that each set of maps compatible with the moduli is also
compatible with \call R.
\end{remark}

\begin{defi} Let $(\tau,\z)$ be an ideal triangulation with moduli of
  $M$. A {\em developing map} for $(\tau,\z)$ is a continuous
 map $D:\wt M\to \M H^3$ such that for each $i$, the map
 $\f_{_{\Delta_i}}$ lifts to a map $\wt\f_{_{\Delta_i}}$ such that
 $D\circ\wt\f_{_{\Delta_i}}$ extends to a map from $\Delta_i$ to
 $\overline{\M H^3}$ which is compatible with $z_i$. Moreover, if
 $\wt\f_{_{\Delta_i}}^1$ and $\wt\f_{_{\Delta_i}}^2$ are two lifts of
 the same map, then we require that there exists an element $\Phi$ of
 Isom$^+(\M H^3)$ such that
$D\circ\wt\f_{_{\Delta_i}}^1= \Phi\circ D\circ\wt\f_{_{\Delta_i}}^2$.
\end{defi}

Let $(\tau,\z)$ be a triangulation with moduli of $M$ and suppose
 \call C holds. Then there exists a representation
 $h:\pi_1(M)\to$Isom$^+(\M H^3)$ which is well-defined up to
 conjugation, defined as follows:

Let $\Delta_{i_1}$ be a base-tetrahedron.
A path of tetrahedra is a sequence
$\gamma=\Delta_{i_1}\stackrel{r_{j_1}}{\dashrightarrow}\Delta_{i_2}\cdots
\stackrel{r_{j_s}}{\dashrightarrow}\Delta_{i_{s+1}}$ 
such that the $\Delta_{i_k}$'s
are tetrahedra of $\tau$ and $r_{j_k}$ is a face pairing rule of
$\tau$  from a 
face of $\Delta_{i_k}$ to one of $\Delta_{i_{k+1}}$.
A loop is a path in which $i_{s+1}=i_1$.
For each $i$ let $g_i:\Delta_i\to\overline{\M H^3}$ be a map
compatible with $z_i$ and for each $j$ let $\psi_j$ the only
orientation preserving isometry which realizes the rule $r_j$ between
the corresponding faces of the geodesic tetrahedra $g_i(\Delta_i)$. We
define $h$ first on the 
set of paths of tetrahedra by setting
$$h(\gamma)=\psi_{j_1}\circ\cdots\circ\psi_{j_s}.$$
Let $P$ be the space of loops of tetrahedra. It is a standard fact
that $P$ project to $\pi_1(M,x_0)$, where $x_0$ is any point of the
image of $\f_{_{\Delta_{i_1}}}$.
The fact that \call C holds imply that $h(\gamma)=h(\gamma')$ if
$[\gamma]=[\gamma']$ in $\pi_1(M,x_0)$. So $h$ induces a map
$h:\pi_1(M,x_0)\to$Isom$^+(\M H^3)$. It is easily checked that $h$ is a
representation and that its conjugacy class does not depend on the
choices of the maps $g_i$ and of the base-tetrahedron.

The representation $h$ is called {\em holonomy} of \z.

\begin{prop}
Let $(\tau,\z)$ be a triangulation with moduli of $M$ such that \z\ is
a solution of equations \call C. Then there exist a developing map
$D$ and a representative of the holonomy such that $D$ is
$\pi_1(M)$-equivariant, where $\pi_1(M)$ acts on \wt M by deck transformations
and on $\M H^3$ via $h$.
\end{prop}

\proof Let \wt\tau\ be the triangulation of \wt M and let
$\Delta_{i_1}$ be a base tetrahedron of \wt\tau. We do the same
construction used to define the holonomy. As above we fix maps $g_i$
(the same $g_i$ for all lifts of the same tetrahedron of $\tau$) and
the isometries $\psi_j$, and as above we define a map \wt h from the space of paths of
tetrahedra of \wt\tau\ to Isom$^+(\M H^3)$.

Now let $\Delta_k$ be a tetrahedron of \wt\tau\ and let $\gamma$ be a
path from $\Delta_{i_1}$ to $\Delta_k$. We define
$D_k:\Delta_k^*\subset\Delta_k\to \M H^3$ by
$D_k=\wt h(\gamma)\circ g_k$. Since \wt M is simply
connected, then  the definition of $D_k$ does not depend
on the choice of $\gamma$. Moreover, since the $g_i$'s are compatible
with the moduli, the union of the $D_k$ projects to a
developing map $D: \wt M\to \M H^3$.

If we fix a base point in \wt M  in the
base-tetrahedron, then a representative $h$ of the holonomy is fixed;
and the $\pi_1(M)$-equivariance of $D$ follows by construction. \qed

\begin{remark}
It is easy to see that if a developing map exists, then
\call C holds and so the holonomy is defined. Moreover,
every developing map can be obtained via the
above construction. It follows that for each developing map $D$ the choices
of a base point and its lift determine a representative $h$ of the
holonomy such that $D$ is $\pi_1(M)$-equivariant.
\end{remark}

\begin{defi}\label{d dm}
Let $(\tau,\z)$ be a triangulation with moduli of $M$ such that \z\ is
a solution of equations \call C.
Let $N$ be an oriented, complete hyperbolic 3-manifold (hence $\wt N=\M H^3$).
We call a map $f:M\to N$ {\em  hyperbolic} if its lift $\wt f:\wt M\to
\M H^3$ is a developing map.
\end{defi}

\begin{remark}\label{r n}
  Let $M, \tau, \z, N, f$ as in Definition~\ref{d dm}. Let $V(\Delta_i)$
  denote the 0-skeleton of $\Delta_i\in\tau$. The maps
  $\f_{_{\Delta_i}}$, and so
  also $f\circ\f_{_{\Delta_i}}$, are defined only on
  $\Delta_i^*$. Nevertheless,
  since \wt f is a developing map, then the maps $\wt
  f\circ\wt\f_{_{\Delta_i}}$
  extend to the whole $\Delta_i$. Thus
  $\wt f\circ\wt\f_{_{\Delta_i}}(V(\Delta_i))$ is well-defined.
\end{remark}

\begin{defi}\label{d spa}
Let $M,\tau,\z,N,f$ be as in Definition~\ref{d dm}.
Let $\gamma$ be an oriented geodesic in N. Let $v$
be a vertex of $\tau$. We say
that $f$ {\em spirals around $\gamma$} near $v$ if, in a suitable
half-space model of $\M H^3$ in which $\wt \gamma$ is the oriented line
$(0,\infty)$, for each tetrahedron $\Delta_i$ having $v$ as a  vertex
there exists a lift $\wt f\circ\wt\f_{_{\Delta_i}}$ which
maps $v$ to $\infty$.
\end{defi}

\begin{prop}\label{p eg}
  Let $M,\tau,\z,N,f$ be as in definition~\ref{d dm}.
Let $\Gamma$ be the subgroup of Isom$(\M H^3)$ such that
$N=\M H^3/\Gamma$.
Then the image of a holonomy representation relative to \wt f
consists of elements of $\Gamma$. Moreover, called $h_N$ the holonomy
of $N$, i.e. the
isomorphism $h_N:\pi_1(N)\to \Gamma$, we have that $h=h_N\circ f_*$.
\end{prop}
\proof
The group $\pi_1(M)$ acts on \wt M via deck transformations and
$\pi_N$ acts on $\wt N=\M H^3$ via $h_N$. 
It is possible to chose the base-points in such a way that for
all $\alpha\in\pi_1(M)$ and for all $x\in\wt M$ we have $\wt
f(\alpha(x))=h_N f_*(\alpha)(\wt f(x))$. Since $\wt f$ is a developing
map, then $\wt f(\alpha(x))=h(\alpha)(\wt f(x))$. It follows that
$h_Nf_*(\alpha)$ and $h(\alpha)$ coincide on the image of \wt f. Since
the dimension of Im$(\wt f)$ is at least two and since both
$h(\alpha)$ and $f_*(\alpha)$ are orientation preserving isometries,
then they coincide on the whole $\M H^3$.\qed

\section{Hyperbolic Dehn filling equations}\label{s hde}

First of all, we recall the definition of Dehn filling of a manifold.

\begin{defi}
  Let M be the interior of an oriented 3-manifold \wb M such that
  $\partial \wb M$ is a union of 
  tori, $\partial \wb M= \sqcup_n T_n$. For each $T_n$ let $(\mu_n,\lambda_n)$
  be a basis for $H_1(T_n,\M Z)$. Let $(p,q)=\{(p_n,q_n)\}$ where
  $(p_n,q_n)$ is either a pair of coprime integers or the symbol
  $\infty$. For each $n$ such that $(p_n,q_n)\neq\infty$, let $L_n$ be
  an oriented solid torus, $m_n$ be a meridian of $T_n'=\partial L_n$,
  $l_n$ be a loop 
  in $T_n$ such that $[l_n]=p_n\mu_n+q_n\lambda_n$ and
  $\f_n:T_n\to T_n'$ be an orientation reversing
  homeomorphism such that
  $\f_n(l_n)=m_i$.

The Dehn filling  of M with parameters $(p,q)$ is the manifold
$$M_{(p,q)}=\wb M\sqcup\{L_n\}\big/\{\f_n\}$$

The tori  $L_n$ are called filling tori.
\end{defi}

\begin{remark}
  The resulting manifold $M_{(p,q)}$ actually depends only on the
  coefficients $(p,q)$ and not on the maps $\f_n$'s.
\end{remark}

\begin{remark}
  Not all the boundary tori are filled in $M_{(p,q)}$. Namely, a torus
    $T_n$ is filled if and only if $(p_n,q_n)\neq\infty$.
If $(p_n,q_n)=\infty$ for all $n$ , then $M_{(p,q)}= M$.
\end{remark}

Suppose that we have an ideal triangulations of a manifold $M$
with a choice of moduli that satisfy \call C. We recall that if
each $z_n$ has positive imaginary part then the moduli define an
 (incomplete) hyperbolic structure on $M$.
 In this section we introduce some  equations on the
moduli, which we call {\em hyperbolic Dehn filling equations}.
When the moduli have positive imaginary part, such equations imply
that the completion of the hyperbolic structure defined by the
moduli on $M$ is a fixed Dehn filling of $M$.

More generally, these equations can be written down even without
restrictions on the imaginary part of the moduli. In that case, in
general, there is not an obvious geometric interpretation of the
solutions of the equations. For this reason we distinguish between
algebraic and geometric solutions of the hyperbolic Dehn filling
equations.

The principal condition expressed by the equations is that if $m$ is a
loop in a boundary torus killed in homology by the filling, then $h(m)=1$.

From now on, let $M$ be the interior of an oriented compact
3-manifold $\wb M$ such that $\partial\wb M$ is the disjoint
union of $k$ tori, let $\tau$ be an ideal triangulation of $M$ and
let $\z$ be a choice of moduli for $\tau$ satisfying the
compatibility equations \call C.

Let $T\subset \partial \wb M$ be a boundary torus.
We push $T$ a little bit in $M$ and we consider
the image of the natural map $\pi_1(T)\to \pi_1(M)$.
We consider the 
half-space model $\M C\times \M R^+$ of $\M H^3$ and a developing
map $D$ such that the vertex relative to $T$ is lifted to a vertex
mapped to $\infty$ via $D$. Then
there exists a choice of the base-points such that the image of the
holonomy
$h(\pi_1(T))$ consists
of maps which fix $\infty$. It follows that by considering the
restriction to $\partial \M H^3\equiv \M C\M P^1$
of the maps in $h(\pi_1(T))$, we
obtain a representation
$h_T$ of $\pi_1(T)$ in the automorphism of \M C. Moreover, since the
restriction to $\partial \M H^3$ of a positive 
isometry is a Moebius transformation, then
$h_T$ actually is a representation $h_T:\pi_1(T)\to $Aff(\M C).
 Since $h$ is well-defined up to
conjugation, then the dilation 
component of $h_T$ is well-defined, which is a representation
$\rho_T:\pi_1(T)\to \M C^*$.

Since $\pi_1(T)$ is Abelian, it follows that $h_T(\pi_1(T))$ consists of
maps which commute with each other. Then it is easy to see that
either they are all translations, or they have a common fixed point.
In the former case we have $\rho_T\equiv 1$.
In the latter case, up to conjugation, we can suppose that the fixed
point is $0$. Thus we get 
$h_T=\rho_T$, in the sense that for all $\alpha\in\pi_1(T)$ and
$\zeta\in\M C$,  $h_T(\alpha)[\zeta]= \rho_T(\alpha)\cdot\zeta$ .
\begin{remark}

In the following, if there are no ambiguities,
by writing
$\rho_T\equiv 1$ we mean that
$h_T(\pi_1(T))$ consists of translations and by $h_T=\rho_T$ we mean
that $h_T(\pi_1(T))$ consists of maps which fix $0$. We notice that
$h,\ h_T$ and  $\rho_T$ depend on $\z$. When we need to emphasize that,
we write $\rho(\z)$ and so on.

\end{remark}

To write the equations, we need to work with $\log(\rho_T)$. In the
following definition we fix a suitable determination of the logarithm
of $\rho_T$. Formally, a boundary torus $T$ of $\wb M$ is not contained
in $M$, so a regular neighborhood $U$ of $T$ is not a sub set of $M$. If
there are no ambiguities, we does not distinguish between $U$ and
$U\cap M$.

\begin{defi}\label{d c} Let $M,\tau,\z$ be as above and let $D$ be
  a developing map. Let $T$ be a boundary component of \wb M and \wt T
  be a lift of T. Consider the  model $\M C\times \M R^+$ of $\M H^3$
  such that the vertex relative to \wt T is mapped to $\infty$.
  Suppose that $h_T=\rho_T$ and suppose that the following condition
  holds:

\begin{itemize}
 \item[] There exists regular neighborhood
 $\wt U\subset \wt M$ of  \wt T
 such that the developed image of $\wt U$
 does not intersect the line $(0,\infty)$.
\end{itemize}

Then we choose a determination of   $\log(\rho_T)$ as follows: let
$H$ be the universal covering of $\M H^3\setminus(0,\infty)$ made by
using the covering $\exp:\M C\to \M C^*$. 

Let $x_0$ and $\wt x_0$ be base points in $T$ and \wt T. Let
$\gamma:[0,1]\to T$ be a loop at $x_0$ and \wt \gamma{} be its lift
starting from $\wt x_0$. After pushing a little $T$ inside $M$, let
$\wt\alpha:[0,1]\to\M C^*$ be the complex component of
$D\circ\wt\gamma$. As $D\circ\wt\gamma$ lifts to $H$, 
the path $\wt\alpha$ lifts to a path $\wb \alpha:[0,1]\to \M C$. 
Since $h_T=\rho_T$, then
$\wt\alpha(1)=\rho_T([\gamma])\cdot\wt\alpha(0)$, and then
$\wb \alpha(1)=\log(\rho_T([\gamma]))+\wb \alpha(0)$. 

The points $\wb\alpha(0)$ and $\wb\alpha(1)$ depend only on the
homotopy class of $\gamma$ and on the
choice of the base-points.
If we change the base-points, then the determination of
$\log(\rho_T([\gamma]))$ changes by a conjugation by translations and
so it is well-defined. 

\end{defi}

Let us now fix a basis $(\mu,\lambda)$ for $H_1(T,\M Z)$ and let $(a,b)$ be
a pair of coprime integers. Consider the Dehn filling of $M$ with
parameters $(a,b)$, i.e. the filling in which the meridian $m$ of the
solid torus is mapped to an oriented loop
homotopic to $a\mu+b\lambda$. So the coefficient $(a,b)$ induces an
orientation of $m$. Since the gluing map inverts the orientations of
the boundary tori, then the core $\gamma$ of the
filling tours is canonically oriented by requiring that $m$ turns
around $\gamma$ by following the right-hand rule in the oriented solid
tours.

We are now ready to give the hyperbolic Dehn filling equations.

\begin{defi}\label{d as}
Let $M$ be the interior of an oriented compact
3-manifold $\wb M$ such that $\partial\wb M$ is the disjoint
union of $k$ tori, let $\tau$ be an ideal triangulation of $M$ and
let $\z$ be a choice of moduli for $\tau$ satisfying the
compatibility equations \call C.

For each boundary torus $T_n$ let
$(\mu_n,\lambda_n)$ be a basis for $H_1(T_n,\M Z)$.
Let $(p,q)=\{(p_n,q_n),\ n=1,\dots,k\}$ be such that
$(p_n,q_n)$ is either a pair of coprime integers or the symbol
$\infty$. Let $\rho_n(\z)$ be the dilation component
of the holonomy of $T_n$ when $\z$ varies on the space of solutions of
the compatibility equations.

We say that $\z$ is an {\em algebraic} solution of the
$(p,q)$-equations if for each $n=1,\dots, k$ we have:
\begin{itemize}
\item[-] If $(p_n,q_n)=\infty$ then $\rho_n(\z)\equiv 1$.
\item[-] If $(p_n,q_n)\neq\infty$ then $h_{T_n}(\z)=\rho_n(\z)$, the condition
  of Definition~\ref{d c}
holds, and

$$
p_n\log(\rho_n(\z)[\mu_n]) + q_n\log(\rho_n(\z)[\lambda_n])=2\pi i.
$$
\end{itemize}

We say that $z$ is a {\em geometric} solution of the
$(p,q)$-equations if, called $N=M_{(p,q)}$ the Dehn filling of $M$
with parameters $(p_n,q_n)$, we have:

\begin{itemize}
\item[a)] $N$ is complete hyperbolic and the cores of filling tori are disjoint
  geodesics $\{\gamma_n\}$.
\item[b)] There exists a proper map $f:M\to N\setminus
  \{\gamma_n\}\subset N$ of degree 1, which is hyperbolic w.r.t. \z.
\item[c)]For each boundary torus $T_n$ with $(p_n,q_n)\neq\infty$,
called $v_n$ the vertex correspondent to $T_n$, $f$ spirals around
  the relative $\gamma_n$ near $v_n$ ($\gamma_n$ has the orientation 
  induced by the Dehn filling coefficient $(p_n,q_n)$).
\end{itemize}
\end{defi}

\begin{remark}
When all the coefficients $(p_n,q_n)$ are $\infty$, then the
system of the $(p,q)$-equations is exactly the classical system
$\call M$ of the so-called {\em completeness equations}. When the
moduli have positive imaginary part the equations \call M imply that
the hyperbolic structure defined by the moduli on $M$ is complete
(of finite volume).
\end{remark}

\begin{teo}\label{t gisa} Let $M,\tau,\{(\mu_n,\lambda_n)\},(p,q)$ be as in
Definition~\ref{d as}. Then
each geometric solution of the $(p,q)$-equations is also algebraic.
\end{teo}

\proof Let \z\ be a geometric solution of the $(p,q)$-equation.

Let $\Gamma$ be the subgroup of Isom$^+(\M H^3)$ such that $N=\M
H^3/\Gamma$. 
So the holonomy of $N$, as a hyperbolic manifold, is an
isomorphism $h_N:\pi_1(N)\to\Gamma$. 
By Proposition~\ref{p eg}
the holonomy of $M$ is obtained by composing
the homomorphism $f_*:\pi_1(M)\to\pi_1(N)$ with $h_N$.

Each geodesic $\gamma_n$ (the cores of the filling tori) can be viewed
as an element of $\pi_1(N)$.
For each $n$, let $\Gamma_n\subset\Gamma$ be the set of all conjugates
of $h_N(\gamma_n)$ and let $P$ be the
set of all parabolic elements of $\Gamma$. 
Note that $P$ is exactly the image of
all boundary elements of $\pi_1(N)$. Since $f$ is
proper, \wt f\ maps each vertex of $\wt\tau$ to a fixed point
of an element either of $\cup_n\Gamma_n$ or of $P$.

Moreover $f$ is surjective on $N\setminus\{\gamma_n\}$ because 
it has degree 1. 
Since $f$ spirals around $\gamma_n$ near
$v_n$, it maps the unfilled cusps of $M$ to the cusps of
$N$. This implies that for $(p_n,q_n)=\infty$ the holonomy of $T_n$
consists of parabolic elements and so $\rho_n\equiv 1$.

If $(p_n,q_n)\neq\infty$ then the image of $h_{T_n}$ is contained in
the subgroup of $\Gamma$ generated by $\gamma_n$, so
$h_{T_n}=\rho_n$. The fact that Im$(f)=N\setminus\{\gamma_n\}$ implies
the condition of Definition~\ref{d c}, and  the fact that $N=M_{(p,q)}$
implies that $p_n\log(\rho_n(\z)[\mu_n]) +
q_n\log(\rho_n(\z)[\lambda_n])=2\pi i$. \qed

For each subgroup $\Gamma$ of Isom$(\M H^3)$, let Fix$(\Gamma)$
denote the set of fixed points of all the elements of $\Gamma$.

\begin{remark}\label{r fgu}
In the proof of Theorem~\ref{t gisa} we showed that
for each $\Delta_n$ of $\tau$, we have $\wt
f\circ\wt\f_{\Delta_n}(V(\Delta_n))\subset$ Fix$(\Gamma)$.
\end{remark}

\begin{remark}
It is well-known that if each $z_n$ has positive imaginary part, then
an algebraic solution is also a geometric solution.
\end{remark}

\begin{remark}
  In Section~\ref{s 5} we give examples of an algebraic solutions
  that are not geometric.
\end{remark}

\section{Uniqueness}\label{s 4}

In this section we prove the uniqueness of the geometric solutions.
We'll need the following version of the rigidity theorem for
complete hyperbolic 3-manifolds of finite volume, which can be found
in~\cite{BCS} and in~\cite{F2}.

\begin{teo}[Strong statement of Mostow's rigidity]\label{t m}
  Let $M_1$ and $M_2$ be two complete connected hyperbolic 3-manifolds
  of finite volume.
  Let $f:M_1\to M_2$ be a continuous proper map such that
  vol$(M_1)=|$deg$(f)|$vol$(M_2)$. Then f is proper homotopic to a locally
  isometric covering of degree deg$(f)$ of $M_1$ onto $M_2$.
\end{teo}

This theorem in particular implies that a 3-manifold carries at most one
hyperbolic structure up to isometries. In the following when we
speak about a hyperbolic 3-manifold $M$, we mean that $M$ is equipped with
its unique hyperbolic structure.

\

In the following we fix an oriented 3-manifold $M$ with $k$ toroidal cusps
equipped with an ideal triangulation $\tau=(\{\Delta_i\},\{r_j\})$ and
a realization ${\cal R}=(\{f_j\},\f)$.
For each boundary torus $T_n$ we fix a basis $(\mu_n,\lambda_n)$ for
$H_1(T_n,\M Z)$.
We fix the Dehn filling coefficients $(p,q)=\{(p_n,q_n)\}$ and a
geometric solution
$\z=\{z_n\}$ of the $(p,q)$-equation.

Let $N=M_{(p,q)}$ equipped with its hyperbolic structure $N=\M
H^3/\Gamma$ with $\Gamma<$Isom$(\M H^3)$ and let $\{\gamma_n\}$ be
the set of the geodesic cores of the filling tori.

The following proposition shows that the set of the geodesics
$\gamma_n$ is uniquely determined by $(p,q)$.

\begin{prop}\label{p 41}
  Let $\mathfrak S_1$ and $\mathfrak S_2$ be two finite-volume,
  complete hyperbolic structure on $N$  such that the $\gamma_n$'s are
  geodesics w.r.t. both  $\mathfrak S_1$ and $\mathfrak S_2$.

Then there exists an isometry $\alpha:(N,\mathfrak S_1)\to(N,\mathfrak
S_2)$ such that for each n $\alpha(\gamma_n)=\gamma_n$.
\end{prop}

\proof By rigidity, the identity $Id:(N,\mathfrak S_1)\to(N,\mathfrak
S_2)$ is homotopic to an isometry $\alpha$. Thus for each $n$ the loop
$\gamma_n$ is homotopic to $\alpha(\gamma_n)$.

By hypothesis $\gamma_n$ is geodesic w.r.t. both $\mathfrak S_1$ and
$\mathfrak S_2$. Since $\alpha$ is an isometry it follows that
$\alpha(\gamma_n)$ is a geodesics w.r.t. $\mathfrak S_2$. Hence
$\gamma_n$ and $\alpha(\gamma_n)$ are geodesics
w.r.t. $\mathfrak S_2$ and are homotopic, and so they must coincide.
\qed

The remaining part of the section is devoted to prove the following

\begin{teo}\label{p uz} In the hypotheses fixed at the beginning of
  this section,
  the moduli $z_i$'s are uniquely determined by the coefficients $(p,q)$.
\end{teo}

\proof
Let $f$ be a hyperbolic map as in Definition~\ref{d as}. By definition
of hyperbolic map, the lift \wt f is a developing map. 
By Proposition~\ref{p eg} the holonomy of \z\ is the composition
$h=h_N\circ f_*$, where $h_N$ is the holonomy of the hyperbolic
manifold $N$. 

Let $v_n$ be the vertex of $\tau$ relative to the
nth cusp of $M$ and let $\wt v_n$ be one of its lifts.  
Let Stab$(\wt v_n)$ be
the stabilizer of $\wt v_n$ in $\pi_1(M)$ 
(acting on \wt M via deck
transformations). 
Then $\wt f(\wt v_n)$ is fixed by $h(\textrm{Stab}(\wt v_n))$.
If $T_n$ is the boundary torus relative to $v_n$, then   
Stab$(\wt v_n)$ is conjugate to $\pi_1(T_n)$. It follows that
$h(\textrm{Stab}(\wt v_n))$ is either parabolic (if $(p_n,q_n)=\infty$)
or generated by a conjugate of $h_N(\gamma_i)$. 
In the former case $h(\textrm{Stab}(\wt v_n))$ consists of exactly one
point and so $\wt f(\wt v_n)$ is completely determined. 
In the latter case $h(\textrm{Stab}(\wt v_n))$ consists of two
points, but the condition $c)$ of Definition~\ref{d as} allows us to
determine $\wt f(\wt v_n)$.
  
Since $\wt f$ is a developing map, then the modulus of a tetrahedron
is completely determined by the $\wt f$-image of its vertices. It
follows that the moduli $z_i$'s depend only on $h_N\circ f_*$. We now prove
that the moduli does not depend on $f$.

Let $g$ be another hyperbolic map as in Definition~\ref{d as}. Let
us call $h_f$ and $h_g$ the holonomy relative respectively to $f$ and
$g$. 
From the properties of the holonomy (see Section~\ref{s itwm}) it
follows that there exists an element $\f\in$ PSL$(2,\M C)$ such that 
$h_g=\f h_f\f^{-1}$.
Thus if a point $p$ is fixed by $h_f(\textrm{Stab}(\wt v_n))$, then
the point $\f(p)$ is fixed by 
$\f h_f(\textrm{Stab}(\wt v_n))\f^{-1}=h_g(\textrm{Stab}(\wt v_n))$.
By using again the condition $c)$ of Definition~\ref{d as}, we see that
$\wt g(v)=\f\wt f(v)$ for each vertex $v$ of \wt\tau. Since $\f$ is an
isometry, then the modulus of a tetrahedron of \wt\tau{} depends only
on $h_N$. In other words the moduli $z_i$'s depend only on the
hyperbolic structure of $N$, which is unique by Theorem~\ref{t m} and
the next Lemma.\qed

\begin{lemma}
The manifold $N$ has finite volume.
\end{lemma}

\proof
Let vol$(\Delta_i)$ denote the hyperbolic volume of the ideal tetrahedron of
modulus $z_i$, where vol$(\Delta_i)$ is taken negative if Im$(z_i)<0$.
Since $f$ is a hyperbolic map, then the volume of its image satisfies
vol$($Im$(f))\leq\sum |$vol$(\Delta_i)|<\infty$.

Moreover $f$ is a degree 1 map from $N\setminus\{\gamma_n\}$ to
$N\setminus\{\gamma_n\}$,
hence vol$(N)=$ vol$(N\setminus\{\gamma_n\})=$vol$($Im$(f))<\infty$.\qed

\section{Examples}\label{s 5}

In this section we do explicit calculations of the solutions of
the compatibility and completeness equation for some particular
one-cusped 3-manifold.

We now fix some notations. Let $L$ and $R$ be the following matrix of
SL$(2,\M Z)$:
$$L=\begin{pmatrix}1&1\\0&1\end{pmatrix}\qquad
R=\begin{pmatrix}1&0\\1&1\end{pmatrix}$$
Each element of $A$ SL$(2,\M Z)$ can be written as a product $A=\prod
A_i^{n_i}$, with $A_i\in\{L,R\}$ and $n_i\in\M N$.  

Let $S$ be the punctured torus $(\M R^2\setminus\M Z^2)/\M Z^2$. 
Then each element $A\in$ SL$(2,\M Z)$
induce an homeomorphism $\f_A$ of $S$. Given $A=\prod A_i^{n_i}$, we
call $\prod A_i^{n_i}$ the manifold obtained from $S\times[0,1]$ by
gluing $(x,0)$ to $(\f_A(x),1)$.

For such manifolds, using the algorithm described in
\cite{fh}, one easily obtain an ideal triangulation with $\sum n_i$ 
tetrahedra. 

We notice that the complement of the figure-eight knot is the manifold
$LR$, and its standard ideal triangulation with two tetrahedra is the
one obtained according with \cite{fh}.

We use the following notations for labeling the simplices.
For each vertex $v$ of a tetrahedron $X$, we call
$X_v$ the triangle obtained by chopping off the vertex $v$ from
$X$ and $X^v$ the face of $X$ opposite to $v$. Given a tetrahedron
$X$ and two vertices $v,w$ of $X$, by abusing notation, we use the
label $w$ also for the edge of the triangle $X_v$ corresponding to
the face $X^w$. A modulus for a tetrahedron $X$ is named $z_X$ and
we will specify the edge to which is referred.

\subsection{The manifold $LR^3$}\label{ss l2r3}

Let $M$ be the manifold $LR^3$, i.e. the manifold obtained as described
above by using the
element
$LR^3=
\begin{pmatrix}
4 & 1
\\ 3 & 1\\
\end{pmatrix}
=
\begin{pmatrix}1&1\\0&1\\
\end{pmatrix}
\begin{pmatrix}1&0\\1&1\\
\end{pmatrix}
\begin{pmatrix}1&0\\1&1\\
\end{pmatrix}
\begin{pmatrix}1&0\\1&1\\
\end{pmatrix}$ of SL$(2,\mathbb Z)$.

 Using the algorithm described in
\cite{fh}, we get the ideal triangulation $\tau$ of $M$ with four
tetrahedra, labeled $A,\,B,\,C,\,D$, pictured in Figure~\ref{f
1}.

\setlength{\unitlength}{1mm}
\begin{figure}[h]
\begin{picture}(150,90)
\put(0,10){
  \put(0,0){\line(0,1){25}}\put(0,13){\qbezier(0,0)(0,0)(-1,-2)\qbezier(0,0)(0,0)(1,-2)}
  \put(0,0){\line(1,0){50}}\multiput(28,0)(-3,0){2}{\qbezier(0,0)(0,0)(-2,-1)\qbezier(0,0)(0,0)(-2,1)}
  \put(50,0){\line(0,1){25}}\put(50,13){\qbezier(0,0)(0,0)(-1,-2)\qbezier(0,0)(0,0)(1,-2)}
  \put(0,25){\line(1,0){50}}\multiput(28,25)(-3,0){2}{\qbezier(0,0)(0,0)(-2,-1)\qbezier(0,0)(0,0)(-2,1)}
  \put(0,0){\line(2,1){50}}\multiput(28,14)(-2.2,-1.1){2}{\qbezier(0,0)(0,0)(-2.5,0)\qbezier(0,0)(0,0)(-1.5,-2)}
  \put(50,0){\line(-2,1){21}}\put(30,10){\qbezier(0,0)(0,0)(2.5,0)\qbezier(0,0)(0,0)(1.5,-2)\qbezier(3.7,0)(3.7,0)(2.7,-2.6)}
  \put(0,25){\line(2,-1){20}}
  \put(-3,-3){$0$}
  \put(51,-3){$\frac{0}{1}$}
  \put(51,26){$\frac{1}{1}$}
  \put(-3,26){$\frac{1}{0}$}

\put(90,0){
  \put(0,0){\line(0,1){25}}\put(0,13){\qbezier(0,0)(0,0)(-1,-2)\qbezier(0,0)(0,0)(1,-2)}
  \put(0,0){\line(1,0){50}}\multiput(28,0)(-3,0){2}{\qbezier(0,0)(0,0)(-2,-1)\qbezier(0,0)(0,0)(-2,1)}
  \put(50,0){\line(0,1){25}}\put(50,13){\qbezier(0,0)(0,0)(-1,-2)\qbezier(0,0)(0,0)(1,-2)}
  \put(0,25){\line(1,0){50}}\multiput(28,25)(-3,0){2}{\qbezier(0,0)(0,0)(-2,-1)\qbezier(0,0)(0,0)(-2,1)}
  \put(0,0){\line(2,1){50}}\multiput(28,14)(-2.2,-1.1){3}{\qbezier(0,0)(0,0)(-2.5,0)\qbezier(0,0)(0,0)(-1.5,-2)}
  \put(50,0){\line(-2,1){20}}
  \put(0,25){\line(2,-1){22}}\multiput(20,15)(-2.2,1.1){2}{\qbezier(0,0)(0,0)(-1.5,2)\qbezier(0,0)(0,0)(-2.5,0)}
  \put(-3,-3){$0$}
  \put(51,-3){$\frac{1}{1}$}
  \put(51,26){$\frac{2}{1}$}
  \put(-3,26){$\frac{1}{0}$}
}

\put(0,45){
  \put(0,0){\line(0,1){25}}\multiput(0,15.5)(0,-3){3}{\qbezier(0,0)(0,0)(-1,-2)\qbezier(0,0)(0,0)(1,-2)}
  \put(0,0){\line(1,0){50}}\multiput(28,0)(-3,0){2}{\qbezier(0,0)(0,0)(-2,-1)\qbezier(0,0)(0,0)(-2,1)}
  \put(50,0){\line(0,1){25}}\multiput(50,15.5)(0,-3){3}{\qbezier(0,0)(0,0)(-1,-2)\qbezier(0,0)(0,0)(1,-2)}
  \put(0,25){\line(1,0){50}}\multiput(28,25)(-3,0){2}{\qbezier(0,0)(0,0)(-2,-1)\qbezier(0,0)(0,0)(-2,1)}
  \put(0,0){\line(2,1){50}}\put(28,14){\qbezier(0,0)(0,0)(-2.5,0)\qbezier(0,0)(0,0)(-1.5,-2)\qbezier(-3.7,0)(-3.7,0)(-2.7,-2.6)}
  \put(50,0){\line(-2,1){21}}\put(30,10){\qbezier(0,0)(0,0)(2.5,0)\qbezier(0,0)(0,0)(1.5,-2)}
  \put(0,25){\line(2,-1){20}}
  \put(-3,-3){$0$}
  \put(51,-3){$\frac{1}{1}$}
  \put(51,26){$\frac{3}{2}$}
  \put(-3,26){$\frac{2}{1}$}
}

  \put(90,45){
  \put(0,0){\line(0,1){25}}\put(0,13){\qbezier(0,0)(0,0)(-1,-2)\qbezier(0,0)(0,0)(1,-2)\qbezier(-1,-3)(-1,-3)(1,-3)}
  \put(0,0){\line(1,0){50}}\multiput(28,0)(-3,0){2}{\qbezier(0,0)(0,0)(-2,-1)\qbezier(0,0)(0,0)(-2,1)}
  \put(50,0){\line(0,1){25}}\put(50,13){\qbezier(0,0)(0,0)(-1,-2)\qbezier(0,0)(0,0)(1,-2)\qbezier(-1,-3)(-1,-3)(1,-3)}
  \put(0,25){\line(1,0){50}}\multiput(28,25)(-3,0){2}{\qbezier(0,0)(0,0)(-2,-1)\qbezier(0,0)(0,0)(-2,1)}
  \put(0,0){\line(2,1){50}}\put(26,13){\qbezier(0,0)(0,0)(-2.5,0)\qbezier(0,0)(0,0)(-1.5,-2)}
  \put(50,0){\line(-2,1){21}}\multiput(30,10)(2.2,-1.1){3}{\qbezier(0,0)(0,0)(2.5,0)\qbezier(0,0)(0,0)(1.5,-2)}
  \put(0,25){\line(2,-1){20}}
  \put(-3,-3){$0$}
  \put(51,-3){$\frac{1}{1}$}
  \put(51,26){$\frac{4}{3}$}
  \put(-3,26){$\frac{3}{2}$}
}

}

 \put(25,0){$A$}
 \put(115,0){$B$}
 \put(25,45){$C$}
 \put(115,45){$D$}
\end{picture}
\caption{Ideal triangulation of $M$}\label{f 1}
\end{figure}
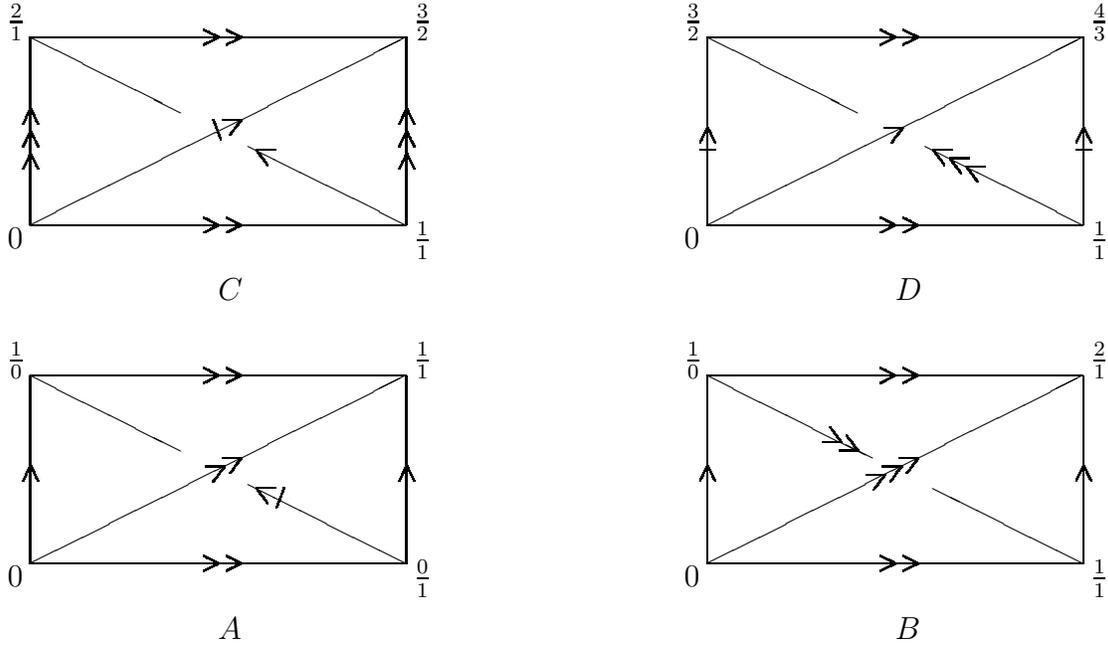

We label the vertices of the tetrahedra as in Figure~\ref{f 1} (we
use such labels because they are natural using the algorithm of
\cite{fh}). The moduli are
 referred to the edge $\overline{0\,\frac{1}{1}}$ (note that
such edge is common to all the tetrahedra).

The face-pairing rules of $\tau$ are, according with the arrows of the
picture:
$$\begin{array}{ccccccc}
 A^{\frac{0}{1}}\longleftrightarrow B^{\frac{2}{1}}&~&
 B^{\frac{1}{0}}\longleftrightarrow C^{\frac{3}{2}}&~&
 C^{\frac{2}{1}}\longleftrightarrow D^{\frac{4}{3}}&~&
 D^{\frac{3}{2}}\longleftrightarrow A^{\frac{1}{1}}\\
 A^{\frac{1}{0}}\longleftrightarrow B^{0}&~&
 B^{\frac{1}{1}}\longleftrightarrow C^{0}&~&
 C^{\frac{1}{1}}\longleftrightarrow D^{0}&~&
 D^{\frac{1}{1}}\longleftrightarrow A^{0}
 \end{array}$$

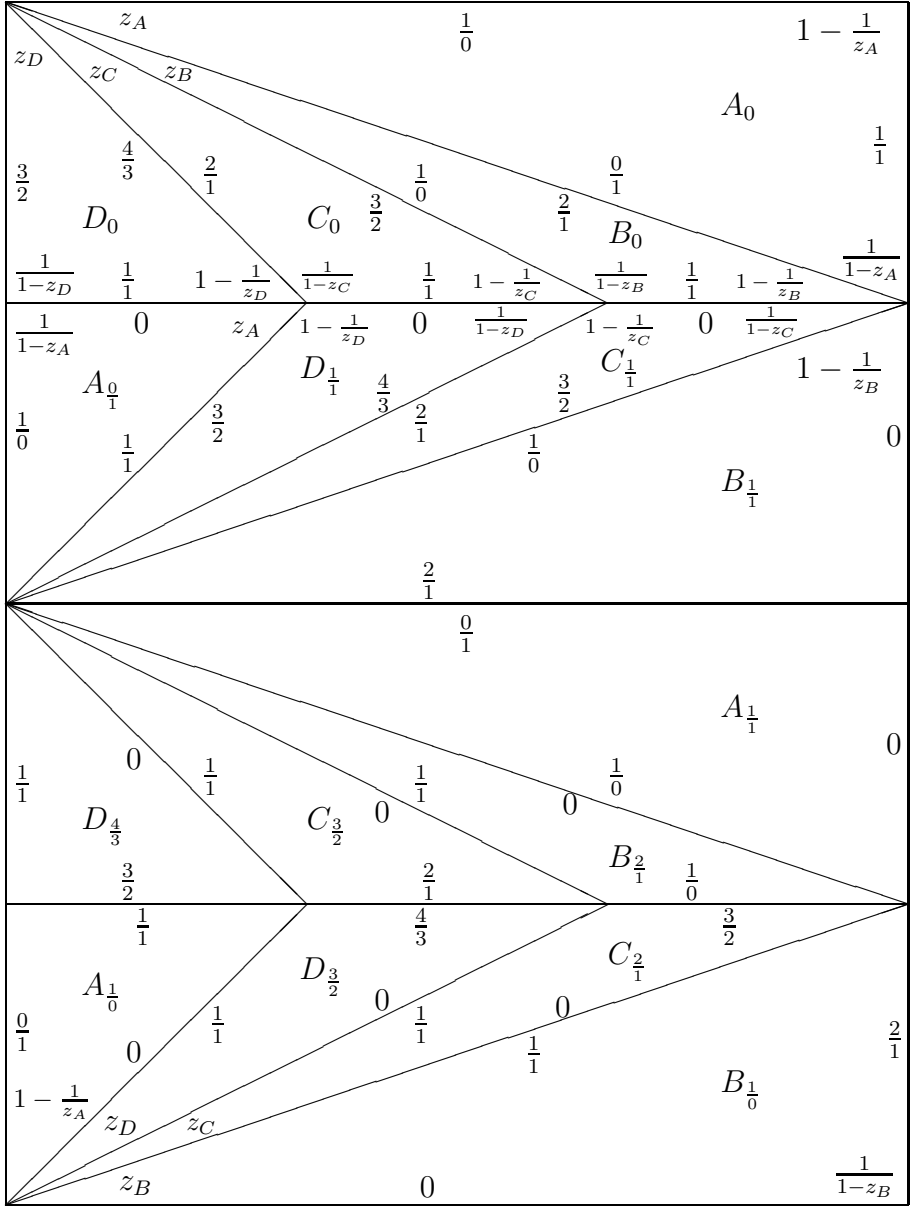
\begin{figure}
\begin{picture}(120,160)
 \multiput(0,0)(0,40){5}{\line(1,0){120}}
 \multiput(0,0)(120,0){2}{\line(0,1){160}}
  \multiput(0,0)(0,80){2}{
    \put(120,40){\line(-3,1){120}}
    \put(0,80){\line(2,-1){80}}
    \put(0,80){\line(1,-1){40}}

    \put(0,0){\line(3,1){120}}
    \put(0,0){\line(2,1){80}}
    \put(0,0){\line(1,1){40}}
  }

 \put(95,145){$A_0$}
 \put(80,128){$B_0$}
 \put(40,130){$C_0$}
 \put(10,130){$D_0$}

     % labels A_0
 \put(115,140){$\frac{1}{1}$}
 \put(60,155){$\frac{1}{0}$}
 \put(80,136){$\frac{0}{1}$}
     % labels B_0
 \put(73,131){$\frac{2}{1}$}
 \put(54,135){$\frac{1}{0}$}
 \put(90,122){$\frac{1}{1}$}
     % labels C_0
 \put(48,131){$\frac{3}{2}$}
 \put(26,136){$\frac{2}{1}$}
 \put(55,122){$\frac{1}{1}$}
     % labels D_0
 \put(15,138){$\frac{4}{3}$}
 \put(1,135){$\frac{3}{2}$}
 \put(15,122){$\frac{1}{1}$}

 \put(79,111){$C_{\frac{1}{1}}$}
 \put(39,110){$D_\frac{1}{1}$}
 \put(10,108){$A_\frac{0}{1}$}
 \put(95,95){$B_\frac{1}{1}$}
     % labels C_11
 \put(92,116){$0$}
 \put(54,103){$\frac{2}{1}$}
 \put(73,107){$\frac{3}{2}$}
     %labels D_11
 \put(54,116){$0$}
 \put(27,103){$\frac{3}{2}$}
 \put(49,107){$\frac{4}{3}$}
     % labels A_01
 \put(17,116){$0$}
 \put(15,99){$\frac{1}{1}$}
 \put(1,102){$\frac{1}{0}$}
     % labels B_11
 \put(69,99){$\frac{1}{0}$}
 \put(117,101){$0$}
 \put(55,82){$\frac{2}{1}$}

\put(0,-80){
 \put(95,145){$A_\frac{1}{1}$}
 \put(80,125){$B_\frac{2}{1}$}
 \put(40,130){$C_\frac{3}{2}$}
 \put(10,130){$D_\frac{4}{3}$}

     % labels A_11
 \put(117,140){$0$}
 \put(60,155){$\frac{0}{1}$}
 \put(80,136){$\frac{1}{0}$}
     % labels B_21
 \put(74,132){$0$}
 \put(54,135){$\frac{1}{1}$}
 \put(90,122){$\frac{1}{0}$}
     % labels C_32
 \put(49,131){$0$}
 \put(26,136){$\frac{1}{1}$}
 \put(55,122){$\frac{2}{1}$}
     % labels D_43
 \put(16,138){$0$}
 \put(1,135){$\frac{1}{1}$}
 \put(15,122){$\frac{3}{2}$}

 \put(80,112){$C_{\frac{2}{1}}$}
 \put(39,110){$D_\frac{3}{2}$}
 \put(10,108){$A_\frac{1}{0}$}
 \put(95,95){$B_\frac{1}{0}$}
     % labels C_21
 \put(95,116){$\frac{3}{2}$}
 \put(54,103){$\frac{1}{1}$}
 \put(73,105){$0$}
     %labels D_32
 \put(54,116){$\frac{4}{3}$}
 \put(27,103){$\frac{1}{1}$}
 \put(49,106){$0$}
     % labels A_10
 \put(17,116){$\frac{1}{1}$}
 \put(16,99){$0$}
 \put(1,102){$\frac{0}{1}$}
     % labels B_10
 \put(69,99){$\frac{1}{1}$}
 \put(117,101){$\frac{2}{1}$}
 \put(55,81){$0$}

}

  % Moduli

\put(15,157){{\footnotesize $z_A$}}
\put(105,155){$1-\frac{1}{z_A}$}
\put(111,125){$\frac{1}{1-z_A}$}
\put(1,115){$\frac{1}{1-z_A}$}
\put(30,116){{\footnotesize $z_A$}}
\put(1,13){{\footnotesize$1-\frac{1}{z_A}$}}

\put(21,150){{\footnotesize $z_B$}}
\put(97,122){{\scriptsize $1-\frac{1}{z_B}$}}
\put(78,123){{\scriptsize $\frac{1}{1-z_B}$}}
\put(105,110){$1-\frac{1}{z_B}$}
\put(15,2){$z_B$}
\put(110,3){$\frac{1}{1-z_B}$}

\put(11,150){{\footnotesize $z_C$}}
\put(62,122){{\scriptsize $1-\frac{1}{z_C}$}}
\put(39,123){{\scriptsize $\frac{1}{1-z_C}$}}
\put(77,116){{\scriptsize $1-\frac{1}{z_C}$}}
\put(98,117){{\scriptsize $\frac{1}{1-z_C}$}}
\put(24,10){{\footnotesize $z_C$}}

\put(1,123){$\frac{1}{1-z_D}$}
\put(25,122){{\footnotesize $1-\frac{1}{z_D}$}}
\put(1,152){{\footnotesize $z_D$}}
\put(39,116){{\scriptsize $1-\frac{1}{z_D}$}}
\put(62,117){{\scriptsize $\frac{1}{1-z_D}$}}
\put(13,10){$z_D$}

\end{picture}
\caption{The triangulation of the boundary torus}\label{f 2}
\end{figure}

The induced triangulation on the boundary torus is the one of
Figure~\ref{f 2}.

We can write down the compatibility and completeness equations. It is
easy to check that ${\cal C+M}$ is equivalent to the following system:

$$
\begin{array}{l}
{\cal{C}}
\left\{
\begin{array}{lcl}
{\cal C }_1.&~& z_A(1-\displaystyle{\frac{1}{z_A})^2 z_D^2 z_C^2
z_B^2\frac{1}{1-z_B}}=1\\\ \\
{\cal C}_2.&~&\displaystyle{
(\frac{1}{1-z_A})^2\frac{1}{1-z_D}(1-\frac{1}{z_B})^2\frac{1}{1-z_C}=1}\\
\ \\
{\cal C}_3.&~&\displaystyle{(1-\frac{1}{z_D})^2\frac{1}{1-z_C}z_A=1}\\
\ \\
{\cal C}_4.&~&\displaystyle{(1-\frac{1}{z_C})^2\frac{1}{1-z_D}\frac{1}{1-z_B}=1}
\end{array}
\right.\\ \ \\
\hskip8mm{\cal M}.\hskip.9cm z_Dz_Cz_B(1-z_A)=1
\end{array}
$$

Moreover, the product of the four equations ${\cal C}$ is exactly
the square of the product of all the moduli, and so it is $1$. So
if three equations are satisfied, then also the remaining one must
be satisfied. It follows that we can discard one of the ${\cal
C}$'s.

We discard ${\cal C}_2$. By using ${\cal M}$ in ${\cal C}_1$ and
then ${\cal C}_1$ in ${\cal C}_4$ and ${\cal M}$ we obtain the
following system of equations, equivalent to ${\cal C+M}$:

\begin{equation}\label{c+m}
\left\{
\begin{array}{lcl}
{\cal M}.&~&\displaystyle{\frac{z_Dz_C(1-z_A)^2}{z_A}}=-1\\ \ \\
{\cal C}_1.&~&\displaystyle{z_A=\frac{1}{1-z_B}}\\ \ \\
{\cal C}_3.&~&\displaystyle{(\frac{z_D-1}{z_D})^2\frac{z_A}{1-z_C}=1}\\
\ \\
{\cal C}_4.&~&\displaystyle{(\frac{z_C-1}{z_C})^2\frac{z_A}{1-z_D}=1}
\end{array}
\right.
\end{equation}

By solving the system, we found four
non-degenerate solutions; one completely positive, giving the
hyperbolic structure of $M$, and other with two negative
tetrahedra, and their conjugates (i.e. the same situations but
with inverted orientations).  The following table contains
numerical approximations of the solutions. Note that even if the
modulus $z_B$ is different from the modulus $z_A$, the equation
${\cal C}_1$ implies that the geometric versions of $A$ and $B$
are isometric. \
\\

\begin{tabular}{|c|c|c|}
\hline
\multicolumn{2}{|c|}{Solution 1}& Volumes\\
\hline
$z_A$&
0.4275047 + i 1.5755666 & 0.9158907\\
\hline
$z_B$&
0.8395957 + i 0.5911691 & 0.9158907\\
\hline
$z_C$&
0.7271548 + i 0.2284421 & 0.5786694\\
\hline
$z_D$&
0.7271548 + i 0.2284421 & 0.5786694\\
\hline
\hline
\multicolumn{2}{|c|}{Solution 2}& Volumes\\
\hline
$z_A$&
1.0724942 + i 0.5921114 & 0.8144270\\
\hline
$z_B$&
0.2854042 + i 0.3945194 & 0.8144270\\
\hline
$z_C$&
-1.7271548 - i 0.6779619 & -0.2398640\\
\hline
$z_D$&
-1.7271548 - i 0.6779619 & -0.2398640\\
\hline
\end{tabular}

\ \\

Note that in the case 2, the total volume is particularly small,
this imply that, even if the identification space is defined, it
cannot be a hyperbolic manifold.

\begin{figure}[h]
\includegraphics[bb = 3.5cm 21.5cm 20cm 25cm]{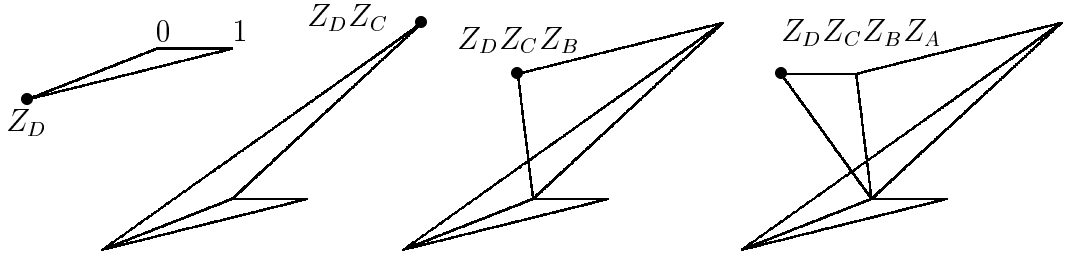}
\caption{The triangles $D_0,\ C_0,\ B_0,\ A_0$
  with the moduli of the solution 2.}\label{f 3}
\end{figure}

\begin{figure}[h]
\includegraphics[bb = 3.7cm 15cm 20cm 20.5cm]{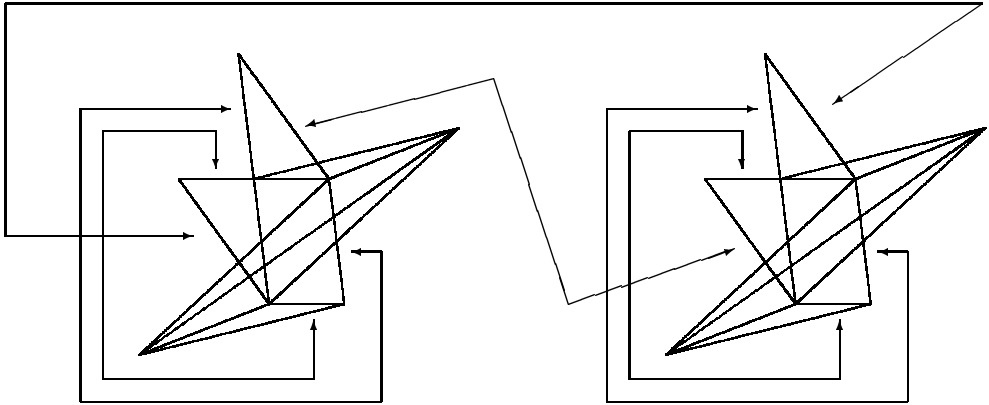}
\caption{Geometric triangulation of the boundary torus, solution 2.}\label{f 4}
\end{figure}

In Figures~\ref{f 3} and~\ref{f 4}, we draw how the triangulations of the
boundary torus of $M$ looks like when we choose the moduli as in
the solution 2. 
There are two types of triangles, the positive ones,
relative to the tetrahedra $A$ and $B$ and the negative ones,
relative to $C$ and $D$. In Figure~\ref{f 3} are pictured the four triangles of
the top strip of the triangulation of Figure~\ref{f 2}.

The two parts of Figure~\ref{f 4} are the top and bottom part of the
triangulation of Figure~\ref{f 2}.

Now we look at the algebraic expression of the solution
of~\ref{c+m}. A simple calculation shows that the moduli can be
expressed by the following equations:

$$
\left\{
\begin{array}{l}
z_C=z_D=w\\ \ \\
z_A=\displaystyle{\frac{w^2}{1-w}}\\ \ \\
z_B=\displaystyle{1-\frac{1}{z_A}}=\displaystyle{\frac{w^2+w-1}{w^2}}\\ \ \\
w^4+2w^3-w^2-3w+2=0
\end{array}
\right.
$$

The four solutions correspond to the four roots $w_1,
\overline{w_1}, w_2, \overline{w_2}$ of the polynomial
$P(w)=w^4+2w^3-w^2-3w+2$. Looking at the reduction (mod 2) of $P$,
we see that $P$ is irreducible over $\mathbb{Z}$, and then also
over $\mathbb{Q}$.

The holonomy representation can be explicitly calculated as a
function of $w$. Let we fix a fundamental domain $F$ for $M$
obtained by taking one copy of each tetrahedron and then
performing the gluing:

$$\begin{array}{ccccccc}
 A^{\frac{1}{0}}\longleftrightarrow B^{0}&~&
 B^{\frac{1}{1}}\longleftrightarrow C^{0}&~&
 C^{\frac{1}{1}}\longleftrightarrow D^{0}
\end{array}$$

Consider  now the geometric version of $F$, i.e. a developed image
of $F$. The holonomy is generated by the isometries corresponding
to the remaining face-pairing rules. We consider the upper
half-space model of $\mathbb{H}^3$ with the coordinates in which
the points $0,1,\infty$ of $\partial \mathbb{H}^3$ are the
vertices of $D$ labeled respectively $\frac{3}{2}, 0,
\frac{4}{3}$. Calculations show that in this model the holonomy is
generated by the elements of $PSL(2,\mathbb{C})$ represented by
the matrices:
$$\begin{pmatrix}1&\frac{w^2}{w^2+w-1}\\0&1\end{pmatrix},~~~~~~~~~
\begin{pmatrix}0&-w\\\frac{1}{w}&-w-1\end{pmatrix},      ~~~~~~~~~
\begin{pmatrix}1&-w^2\\-1&w^2+w-1\end{pmatrix}$$

that respectively correspond to the face-pairing rules
$$\begin{array}{ccccccc}
 A^{0}\longrightarrow D^{\frac{1}{1}}&~&
 C^{\frac{2}{1}}\longrightarrow D^{\frac{4}{3}}&~&
 B^{\frac{2}{1}}\longrightarrow A^{\frac{0}{1}}
 \end{array}$$

What is important is that the entries of such matrices are numbers
belonging to $\mathbb{Q}(w)$ (and this can be proved even without
the explicit calculations).

\begin{prop}\label{p nd}
The solution 2 is not geometric.
\end{prop}

\proof This obviously follows from the uniqueness of geometric solutions,
nevertheless we give an alternative proof.
Let $w_1$ (resp.
$w_2$) be the root of $P$ relative to the solution 1 (resp. 2)
of $\cal{C+M}$. So $w_1$ gives the hyperbolic structure of $M$.
Let $h_j:\pi_1(M)\to PSL(2,\mathbb{C})$ be the holonomy
representation relative to $w_j$, $j=1,2$. Since $P$ is
irreducible and since the entries of the holonomy-matrices are in
$\mathbb{Q}(w)$, it follows that a relation between elements holds
for $h_1$ if and only if it holds for $h_2$. Since $h_1$ is the
holonomy of the complete hyperbolic structure of $M$, then it is
faith-full, and it follows that also $h_2$ is faith-full. 

The image of
$h_2$ cannot be discrete because otherwise $\mathbb{H}^3/h_2$ will
be a hyperbolic manifold $M'$ with a volume too small (actually,
to obtain an absurd it suffices that $vol(h_2)\neq vol(h_1)$).

By Proposition~\ref{p eg} the holonomy of any geometric solution is
discrete, so the solution 2 cannot be geometric.\qed

Similarly, from the fact that $h_2$ is not discrete and
Proposition~\ref{p eg} it follows that there is 
no map, which is hyperbolic w.r.t. the solution 2,
from $LR^3$ to any hyperbolic manifold.

Finally, we show that the image of $h_2$ is dense in $PSL(2,\M C)$. 
We'll need the following standard fact about $PSL(2,\M C)$ (see for
example~\cite{k} or~\cite{G}). 

\begin{lemma}\label{l dns}
Let G be a non-elementary subgroup of $PSL(2,\M C)$ and suppose
that G is not discrete. Then the closure of G is either $PSL(2,\M C)$ or
it is conjugate to $PSL(2,\M R)$ or to a $\M Z_2$-extension of
$PSL(2,\M R)$.
\end{lemma}

\begin{prop}
The image of the holonomy relative to the solution 2 is dense in
$PSL(2,\mathbb{C})$.
\end{prop}

\proof 
It is easy to check that the image of $h_2$ is a
non-elementary subgroup of $PSL(2,\mathbb{C})$.
Suppose that its closure is conjugate to $PSL(2,\M R)$ or to a $\M
Z_2$-extension of $PSL(2,\M R)$. Then there exist a line in $\M
C\cup\{\infty\}=\partial \M H^3$ which is $h_2$-invariant.
By looking at the parabolic elements of $h_2$, it is easy to see that
such a line does not exist. The thesis follows by
Lemma~\ref{l dns}.\qed

This example is interesting for several reasons. On one hand
 it shows that an algebraic solution of $\cal{C+M}$ can be non
 geometric. On the other hand it shows that there is no uniqueness of
 the algebraic solutions. 

Moreover this example does not involve flat tetrahedra, so it is quite
``regular''. Finally, the bad solution of $\call{C+M}$ of $LR^3$ has
the property that  
``everything works OK at the boundary'', in the meaning that the
triangulation with moduli induced on the boundary torus defines on it a
Euclidean structure (up to scaling) with non-zero area. Roughly
speaking, this means that the cusp of $LR^3$ would like to have a
complete hyperbolic structure of finite volume according to the bad
solution of $\call{C+M}$, but the rest of the manifold does not agree.

\subsection{The manifold $L^2R^3$}

The case of $LR^3$ is not at all isolated in the family of fiber
bundle with fiber a punctured torus. In this section we do
calculations for the manifold $L^2R^3$. 

$$\textrm{L}^2\textrm{R}^3
 =
\begin{pmatrix}1&1\\0&1\\
\end{pmatrix}
\begin{pmatrix}1&1\\0&1\\
\end{pmatrix}
\begin{pmatrix}1&0\\1&1\\
\end{pmatrix}
\begin{pmatrix}1&0\\1&1\\
\end{pmatrix}
\begin{pmatrix}1&0\\1&1\\
\end{pmatrix}=\begin{pmatrix}
7 & 2
\\ 3 & 1\\
\end{pmatrix}.$$

Using the algorithm described in
\cite{fh}, we get the ideal triangulation $\tau$ of $M$ with five
tetrahedra, labeled $A,\,B,\,C,\,D,\,E$, pictured in 
Figure~\ref{f 521}. 

\begin{center}
\setlength{\unitlength}{0.6mm}
\begin{figure}[h]\makebox[\textwidth]{
\begin{picture}(200,95)\footnotesize

\put(0,20){

%%%%%%%%%%%%%%% AAAAAAAAAAAAAAAAAAAAAA
\put(0,6){
\put(0,0){\line(0,1){25}}\put(0,13){\qbezier(0,0)(0,0)(-1,-2)\qbezier(0,0)(0,0)(1,-2)}
  \put(0,0){\line(1,0){50}}\multiput(28,0)(-3,0){2}{\qbezier(0,0)(0,0)(-2,-1)\qbezier(0,0)(0,0)(-2,1)}
  \put(50,0){\line(0,1){25}}\put(50,13){\qbezier(0,0)(0,0)(-1,-2)\qbezier(0,0)(0,0)(1,-2)}
  \put(0,25){\line(1,0){50}}\multiput(28,25)(-3,0){2}{\qbezier(0,0)(0,0)(-2,-1)\qbezier(0,0)(0,0)(-2,1)}
  \put(0,0){\line(2,1){50}}\multiput(28,14)(-2.2,-1.1){2}{\qbezier(0,0)(0,0)(-2.5,0)\qbezier(0,0)(0,0)(-1.5,-2)}\put(21,9.6){$\bullet$}
  \put(50,0){\line(-2,1){21}}\put(30,10){\qbezier(0,0)(0,0)(2.5,0)\qbezier(0,0)(0,0)(1.5,-2)\qbezier(3.7,0)(3.7,0)(2.7,-2.6)}
  \put(0,25){\line(2,-1){20}}
  \put(-3,-3){$0$}
  \put(51,-3){$\frac{0}{1}$}
  \put(51,26){$\frac{1}{1}$}
  \put(-3,26){$\frac{1}{0}$}

 }

%%%%%%%%%%%%%%%%% BBBBBBBBBBBBBBBBBBBBBBBBBBB
\put(0,47){

\put(0,0){\line(0,1){25}}\put(0,13){\qbezier(0,0)(0,0)(-1,-2)\qbezier(0,0)(0,0)(1,-2)}
  \put(0,0){\line(1,0){50}}\multiput(28,0)(-3,0){2}{\qbezier(0,0)(0,0)(-2,-1)\qbezier(0,0)(0,0)(-2,1)}\put(20,-1.7){$\bullet$}
  \put(50,0){\line(0,1){25}}\put(50,13){\qbezier(0,0)(0,0)(-1,-2)\qbezier(0,0)(0,0)(1,-2)}
  \put(0,25){\line(1,0){50}}\multiput(28,25)(-3,0){2}{\qbezier(0,0)(0,0)(-2,-1)\qbezier(0,0)(0,0)(-2,1)}\put(20,23.5){$\bullet$}
  \put(0,0){\line(2,1){50}}\multiput(28,14)(-2.2,-1.1){2}{\qbezier(0,0)(0,0)(-2.5,0)\qbezier(0,0)(0,0)(-1.5,-2)}
  \put(50,0){\line(-2,1){21}}\multiput(34,8)(-2.2,1.1){2}{\qbezier(0,0)(0,0)(-1.5,2)\qbezier(0,0)(0,0)(-2.5,0)}
  \put(0,25){\line(2,-1){20}}
  \put(-3,-3){$0$}
  \put(51,-3){$\frac{1}{1}$}
  \put(51,26){$\frac{2}{1}$}
  \put(-3,26){$\frac{1}{0}$}

%%%%%%%%%%%%%%%%%%%%% CCCCCCCCCCCCCCCCCCCCCCCCCCCCC
\put(80,-41){
  \put(0,0){\line(0,1){25}}\put(0,13){\qbezier(0,0)(0,0)(-1,-2)\qbezier(0,0)(0,0)(1,-2)}
  \put(0,0){\line(1,0){50}}\multiput(28,0)(-3,0){2}{\qbezier(0,0)(0,0)(-2,-1)\qbezier(0,0)(0,0)(-2,1)}
  \put(50,0){\line(0,1){25}}\put(50,13){\qbezier(0,0)(0,0)(-1,-2)\qbezier(0,0)(0,0)(1,-2)}
  \put(0,25){\line(1,0){50}}\multiput(28,25)(-3,0){2}{\qbezier(0,0)(0,0)(-2,-1)\qbezier(0,0)(0,0)(-2,1)}
  \put(0,0){\line(2,1){50}}\multiput(28,14)(-2.2,-1.1){3}{\qbezier(0,0)(0,0)(-2.5,0)\qbezier(0,0)(0,0)(-1.5,-2)}
  \put(50,0){\line(-2,1){20}}
  \put(0,25){\line(2,-1){22}}\multiput(20,15)(-2.2,1.1){2}{\qbezier(0,0)(0,0)(-1.5,2)\qbezier(0,0)(0,0)(-2.5,0)}\put(13,16.2){$\bullet$}
  \put(-3,-3){$0$}
  \put(51,-3){$\frac{2}{1}$}
  \put(51,26){$\frac{3}{1}$}
  \put(-3,26){$\frac{1}{0}$}
}
%%%%%%%%%%%%%%%%%%%%%%%%%%%%%%%%%%%%%%%%%%%%%%%%%%%%%%%%%%%%%%%%%%%%%%%%%%%

%%%%%%%%%%%%%%%%%%%%%%%%%% DDDDDDDDDDDDDDDDDDDDDDDDD
\put(80,0){
  \put(0,0){\line(0,1){25}}\multiput(0,15.5)(0,-3){3}{\qbezier(0,0)(0,0)(-1,-2)\qbezier(0,0)(0,0)(1,-2)}
  \put(0,0){\line(1,0){50}}\multiput(28,0)(-3,0){2}{\qbezier(0,0)(0,0)(-2,-1)\qbezier(0,0)(0,0)(-2,1)}
  \put(50,0){\line(0,1){25}}\multiput(50,15.5)(0,-3){3}{\qbezier(0,0)(0,0)(-1,-2)\qbezier(0,0)(0,0)(1,-2)}
  \put(0,25){\line(1,0){50}}\multiput(28,25)(-3,0){2}{\qbezier(0,0)(0,0)(-2,-1)\qbezier(0,0)(0,0)(-2,1)}
  \put(0,0){\line(2,1){50}}\put(28,14){\qbezier(0,0)(0,0)(-2.5,0)\qbezier(0,0)(0,0)(-1.5,-2)\qbezier(-3.7,0)(-3.7,0)(-2.7,-2.6)}
  \put(50,0){\line(-2,1){21}}\put(30,10){\qbezier(0,0)(0,0)(2.5,0)\qbezier(0,0)(0,0)(1.5,-2)}
  \put(0,25){\line(2,-1){20}}
  \put(-3,-3){$0$}
  \put(51,-3){$\frac{2}{1}$}
  \put(51,26){$\frac{5}{2}$}
  \put(-3,26){$\frac{3}{1}$}
}

%%%%%%%%%%%%%%%%%%%%%%%%% EEEEEEEEEEEEEEEEEEEEEEEEEEEEEEE
  \put(160,0){
  \put(0,0){\line(0,1){25}}\put(0,13){\qbezier(0,0)(0,0)(-1,-2)\qbezier(0,0)(0,0)(1,-2)\qbezier(-1,-3)(-1,-3)(1,-3)}
  \put(0,0){\line(1,0){50}}\multiput(28,0)(-3,0){2}{\qbezier(0,0)(0,0)(-2,-1)\qbezier(0,0)(0,0)(-2,1)}
  \put(50,0){\line(0,1){25}}\put(50,13){\qbezier(0,0)(0,0)(-1,-2)\qbezier(0,0)(0,0)(1,-2)\qbezier(-1,-3)(-1,-3)(1,-3)}
  \put(0,25){\line(1,0){50}}\multiput(28,25)(-3,0){2}{\qbezier(0,0)(0,0)(-2,-1)\qbezier(0,0)(0,0)(-2,1)}
  \put(0,0){\line(2,1){50}}\put(26,13){\qbezier(0,0)(0,0)(-2.5,0)\qbezier(0,0)(0,0)(-1.5,-2)}
  \put(50,0){\line(-2,1){21}}\multiput(30,10)(2.2,-1.1){3}{\qbezier(0,0)(0,0)(2.5,0)\qbezier(0,0)(0,0)(1.5,-2)}
  \put(0,25){\line(2,-1){20}}
  \put(-3,-3){$0$}
  \put(51,-3){$\frac{2}{1}$}
  \put(51,26){$\frac{7}{3}$}
  \put(-3,26){$\frac{5}{2}$}
}
}

\put(0,43){
 \put(25,-4){$B$}
 \put(103,-45){$C$}
 \put(103,-4){$D$}
 \put(183,-4){$E$}
\put(25,-45){$A$}
 }

% questa graffa va col put 0 10 iniziale
}
\put(18,40){\parbox{20cm}{\center Face-pairing rules\\ (according with
    the rows)}}
\put(164,30){\line(1,0){40}}
\put(185,30){\vector(-2,-1){15}}

\put(-65,6){
  \parbox{20cm}{\center
        $$\begin{array}{ccccccccc}
        A^{\frac{0}{1}}\leftrightarrow B^{\frac{2}{1}}&~&
        B^{\frac{1}{1}}\leftrightarrow C^{\frac{3}{1}}&~&
        C^{\frac{1}{0}}\leftrightarrow D^{\frac{5}{2}}&~&
        D^{\frac{3}{1}}\leftrightarrow E^{\frac{7}{3}}&~&
        E^{\frac{5}{2}}\leftrightarrow A^{\frac{1}{1}}\\
        A^{\frac{1}{0}}\leftrightarrow B^{0}&~&
        B^{\frac{1}{0}}\leftrightarrow C^{0}&~&
        C^{\frac{2}{1}}\leftrightarrow D^{0}&~&
        D^{\frac{2}{1}}\leftrightarrow E^{0}&~&
        E^{\frac{2}{1}}\leftrightarrow A^{0}
      \end{array}$$
    }
  }

\end{picture}
}
\caption{Ideal triangulation of $M$}\label{f 521}
\end{figure}
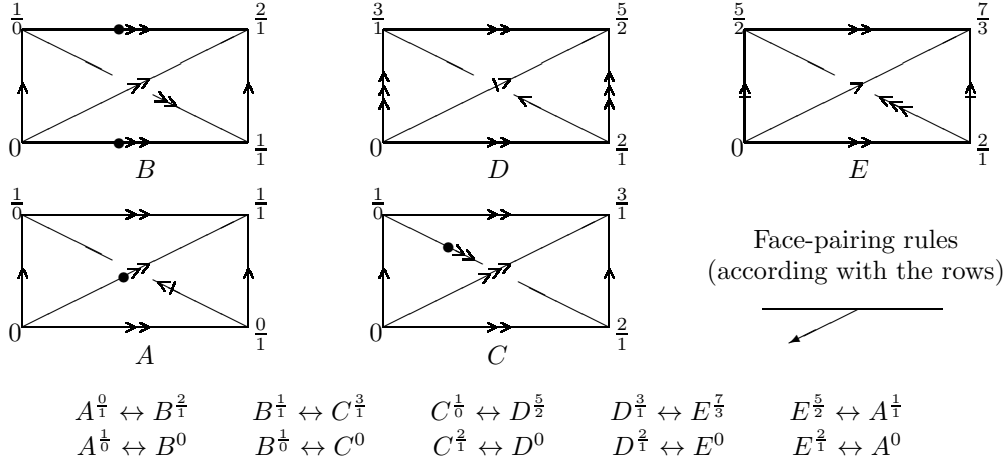
\end{center}

We label the vertices of the tetrahedra as in 
Figure~\ref{f 521}.
The moduli $z_A$
and $z_B$ are referred to the edge $\overline{0\,\frac{1}{0}}$; while
$z_C,\,z_D,\,z_E$ to the edge $\overline{0\,\frac{2}{1}}$.

\begin{figure}[ht]
  \centerline{
\includegraphics[bb = 100 350 500 720]{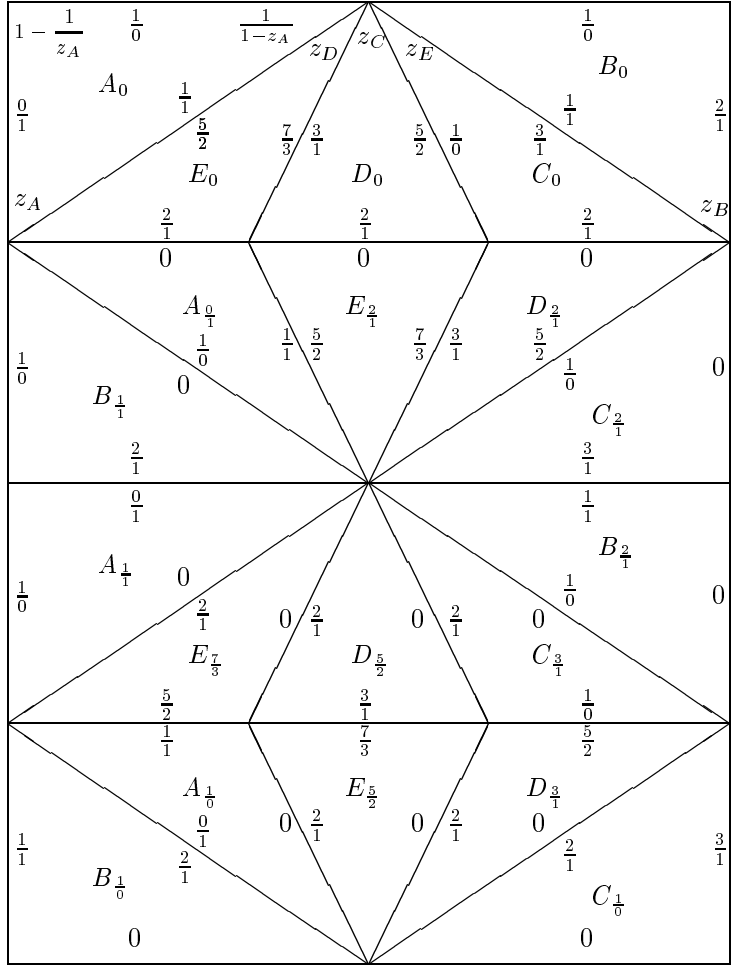}
  }
\caption{The triangulation of the boundary torus}\label{f 522}
\end{figure}

The induced triangulation on the boundary torus is the one of
Figure~\ref{f 522}.

It is easy to see that the system of compatibility and completeness
equations ${\cal C+M}$ is equivalent to the following one: 

\begin{displaymath}
\left\{
\begin{array}{l}
z_Az_B=z_Cz_Dz_E\\
z_C=\displaystyle{\frac{1}{1-z_A}}\\
(1-z_D)^2z_E^2=(1-z_E)^2z_D^2\\
\displaystyle{
(1-\frac{1}{z_A})^2=(1-z_B)^2}\\
\displaystyle{(1-\frac{1}{z_E})^2\frac{1}{1-z_D}(1-\frac{1}{z_A})=1}

\end{array}
\right.
\end{displaymath}

By solving such a system (we have done that with a computer), we found eight
solutions.
The following tables contain numerical approximations of the
solutions. Note that even if the 
modulus $z_A$ is different from the modulus $z_C$, the second equation
implies that the geometric versions of $A$ and $C$
are isometric.
\\
{\footnotesize
\begin{center}
\begin{tabular}{|c|c|c|c|c|}
\cline{3-3}\cline{5-5}
\hline
\multicolumn{2}{|c|}{Solutions 1}& volume &

\multicolumn{1}{||c|}{ Solutions 2}&volume\\
\hline
$z_A$&
$0.75 + i 0.6614378 $& 0.9626730&
\multicolumn{1}{||c|}{$0.75 - i 0.6614378$}& -0.9626730\\
\hline
$z_B$&
$1.25 + i 0.6614378 $& 0.7413987&
\multicolumn{1}{||c|}{$1.25 - i 0.6614378$}& -0.7413987\\
\hline
$z_C$&
$0.5 + i 1.3228756 $& 0.9626730&
\multicolumn{1}{||c|}{$0.5 - i 1.3228756$}& -0.9626730\\
\hline
$z_D$&
$1 $& $*$&
\multicolumn{1}{||c|}{$1$}& $*$\\
\hline
$z_E$&
$1 $& $*$&
\multicolumn{1}{||c|}{$1$}&$*$\\
\hline
\end{tabular}
\end{center}

\begin{center}
\begin{tabular}{|c|c|c|c|c|}
\cline{3-3}\cline{5-5}
\hline
\multicolumn{2}{|c|}{Solutions 3}& volume & 
\multicolumn{1}{||c|}{ Solutions 4}& volume\\
\hline
$z_A$&
$1.588633261 $& 0&
\multicolumn{1}{||c|}{$1.127804076$}&0\\
\hline
$z_B$&
$1.370528159 $&0&
\multicolumn{1}{||c|}{$1.113321168$}& 0\\
\hline
$z_C$&
$-1.69885025 $& 0&
\multicolumn{1}{||c|}{$-7.824476637$}& 0\\
\hline
$z_D$&
$0.3783840018 $& 0&
\multicolumn{1}{||c|}{$0.2518509745$}& 0\\
\hline
$z_E$&
$-3.387066549 $&0&
\multicolumn{1}{||c|}{$-0.6371698130$}& 0\\
\hline
\end{tabular}
\end{center}

\begin{center}
\begin{tabular}{|c|c|c|c|c|}
\cline{3-3}\cline{5-5}
\hline
\multicolumn{2}{|c|}{Solutions 5}& volume& 
\multicolumn{1}{||c|}{ Solutions 6}& volume\\
\hline
$z_A$&
$0.4950484+i0.3298695 $& 0.7399514&
\multicolumn{1}{||c|}{$0.4950484-i0.3298695$}& -0.7399514\\
\hline
$z_B$&
$0.6011109+i0.9321327$& 1.0089809&
\multicolumn{1}{||c|}{$0.6011109-i0.9321327$}& -1.0089809\\
\hline
$z_C$&
$1.3880304+i0.9067580$& 0.7399514&
\multicolumn{1}{||c|}{$1.3880304-i0.9067580$}&-0.7399514\\
\hline
$z_D$&
$0.5022247+i0.2691269$& 0.6433681&
\multicolumn{1}{||c|}{$0.5022247-i0.2691269$}& -0.6433681\\
\hline
$z_E$&
$0.6077815+i0.3441339$& 0.7596486&
\multicolumn{1}{||c|}{$0.6077815-i0.3441339$}& -0.7596486\\
\hline
\end{tabular}
\end{center}

\begin{center}
\begin{tabular}{|c|c|c|c|c|}
\cline{3-3}\cline{5-5}
\hline
\multicolumn{2}{|c|}{Solutions 7}& volume& 
\multicolumn{1}{||c|}{ Solutions 8}&volume\\
\hline
$z_A$&
$0.1467328+i1.2472524$& 0.9386051&
\multicolumn{1}{||c|}{$0.1467328-i1.2472524$}& -0.9386051\\
\hline
$z_B$&
$1.9069644+i0.7908171 $& 0.4782906&
\multicolumn{1}{||c|}{$1.9069644-i0.7908171$}& -0.4782906\\
\hline
$z_C$&
$0.3736330+i0.5461534 $& 0.9386051&
\multicolumn{1}{||c|}{$0.3736330-i0.5461534 $}& -0.9386051\\
\hline
$z_D$&
$1.1826577-i2.5849142$&-0.7155138&
\multicolumn{1}{||c|}{$1.1826577+i2.5849142$}&0.7155138\\
\hline
$z_E$&
$-0.5956636+i1.2429350$& 0.7019645&
\multicolumn{1}{||c|}{$-0.5956636-i1.2429350$}& -0.7019645\\
\hline

\end{tabular}
\end{center}
}
\ \\

The Solutions 1 and 2 contain degenerated tetrahedra. We notice that
the non-degenerate moduli of such solutions are exactly the ones that give the
hyperbolic structure on the manifold obtained by removing the
tetrahedra $D$ and $E$ and adding the gluing rules:
$$\begin{array}{lcl}
C^{\frac{1}{0}}\leftrightarrow A^{\frac{1}{1}}&\textrm{via}&
(0,\frac{3}{1},\frac{2}{1})\leftrightarrow(0,\frac{1}{0},\frac{0}{1})\\
C^{\frac{2}{1}}\leftrightarrow A^{0}&\textrm{via}&
(0,\frac{1}{0},\frac{3}{1})\leftrightarrow(\frac{0}{1},\frac{1}{0},\frac{1}{1}).
\end{array}$$ 

Now we look at the algebraic expression of the Solutions 3-8. 
Let $P(x)=x^6+4x^5+3x^4+3x^3-4x^2+2$.
A simple calculation shows that the moduli can be
expressed in terms of roots of $P$ by the following expressions:

$$
\left\{
\begin{array}{l}
z_A=\displaystyle{\frac{1}{22}}(5w^5+19w^4+9w^3+6w^2-8w+17)\\
z_B=\displaystyle{\frac{1}{44}}(10w^5+49w^4+62w^3+34w^2-16w+34)\\
z_C=\displaystyle{\frac{1}{11}}(-12w^5-39w^4-4w^3-10w^2+72w-32)\\
z_D=\displaystyle{\frac{1}{22}}(-4w^5-13w^4+6w^3+15w^2+2w+4)\\
z_E=w\\
P(w)=0
\end{array}
\right.
$$

As in the case of $LR^3$ the polynomial $P$ is irreducible, and the
solutions $3,4,7,8$ are not geometric.

\subsection{A manifold with non-trivial JSJ decomposition}
The manifold we consider in this section is obtained by gluing to
boundary torus of the complement of the
figure-eight knot a Seifert manifold. 
The resulting manifold, which we call $M$, clearly is not hyperbolic
because it contains an incompressible tours (the old boundary torus).

This example is interesting because the  manifold $M$ admits 
an ideal triangulation with four tetrahedra
such that there exists a positive, partially flat solution of $\call
C+\call M$. Obviously such a solution cannot be geometric, as $M$ is
not hyperbolic. 

At a first look, this seems
 to contradict the result of~\cite{PW}, but fortunately there are no
 contradictions. 
 The point is that
in the present example the moduli 
do not satisfy the equations on the angles. Namely,
when a modulus is positive then it is well defined the angle
associated to a modulus, in such a way that the sum of angles of any
horospherical triangle is always $2\pi$. Then in addiction to
equations $\call C$ one can require that the sum of the angles around
any edge is exactly $2\pi$. Such equations are called $\call C^*$.

In~\cite{PW} is proved that a partially flat solution of $\call
C^*+\call M$ leads to a hyperbolic structure. Here we show a partially
flat solution of $\call C+\call M$ that does not satisfy $\call C^*$. 

This shows that the equations $\call C^*$ play a fundamental role in
order to have hyperbolicity. 
Nevertheless, we notice that equations $\call C^*$ are
not necessary. Namely, there exist examples of
ideal triangulations of the complement of the figure-eight knot whose
unique geometric solution does not satisfy $\call C^*$.

We describe now our manifold $M$.
Let $A$ be the following subset of $\M C$:
$$A=\{z\in\M C: |z|\leq 4,\ |z-2|>1, |z+2|>1\}.$$
A is a disc with two holes. Let $I\subset A$ be the set of the point
with zero real part. 
Let $S$ be the space
obtained from $A\times[0,1]$ by gluing $(z,0)$ to $(-z,1)$ and let
$L$
be the M\"obius strip coming from $I$.
This is the piece we want to glue to complement of 
the figure-eight knot. We call $C_e$ and $C_i$ the external and
internal components of $\partial S$. Note that $\partial L\subset C_e$.   

We will glue $C_e$ to the boundary torus of the complement of the
figure-height knot. To do this, we specify where we glue the boundary of
the M\"obius strip.
Using the classical triangulation of the complement
of the figure-eight knot we get the following picture when looking
from the cusp:
\setlength{\unitlength}{4pt}
\begin{center}
\begin{figure}[ht]
  \begin{picture}(90,18)
    \put(0,0){\line(1,0){80}}
    \put(10,17){\line(1,0){80}}
    \multiput(0,0)(20,0){5}{\qbezier(0,0)(0,0)(10,17)}
    \multiput(20,0)(20,0){4}{\qbezier(0,0)(0,0)(-10,17)}

    \multiput(10,6)(20,0){4}{\pbez{10}(0,0)(-2.5,1.25)(-5,2.5)}
    \multiput(10,6)(20,0){4}{\pbez{10}(0,0)(2.5,1.25)(5,2.5)}
    \multiput(10,6)(20,0){4}{\pbez{10}(0,0)(0,-3)(0,-6)}

    \multiput(20,11)(20,0){4}{\pbez10(0,0)(-2.5,-1.25)(-5,-2.5)}
    \multiput(20,11)(20,0){4}{\pbez{10}(0,0)(2.5,-1.25)(5,-2.5)}
    \multiput(20,11)(20,0){4}{\pbez{10}(0,0)(0,3)(0,6)}

    \setpenn{{\large$.$}}
\put(5,0){\pbez{200}(0,0)(5,8.5)(10,17)}

  \end{picture}
  \caption{The boundary of the complement of the figure-eight knot}
  \label{f bcb}
\end{figure}
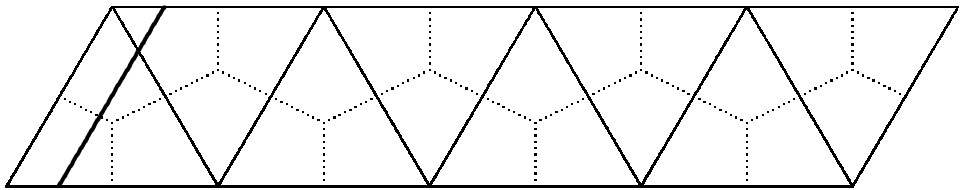
\end{center}

There are pictured the eight equilateral triangles of the
boundary. The dash lines represent the standard spine dual of the ideal
triangulation. Finally the marked line is where we glue $\partial L$.

Since $S$ retracts to $C_e\cup L$, then a spine of $M$ is obtained
simply by gluing a M\"obius strip to the spine of the complement of
the figure-eight knot as in Figure~\ref{f bcb}. Such a spine as a
vertex more than the old one, but is not
standard. 

We perform a {\em lune} move along the M\"obius
strip and  we obtain a standard spine (see~\cite{M} for
details about the moves on the spines) with five vertices. 
As the new spine is standard,
its dual is an ideal triangulation with five tetrahedra. 
Such a triangulation can be simplified with a {\em MP} move, replacing
the three new tetrahedra with an equivalent pair of tetrahedra.

At the end, we get the triangulation of $M$ sketched in 
Figure~\ref{f taucb}.
 
\setlength{\unitlength}{4pt}
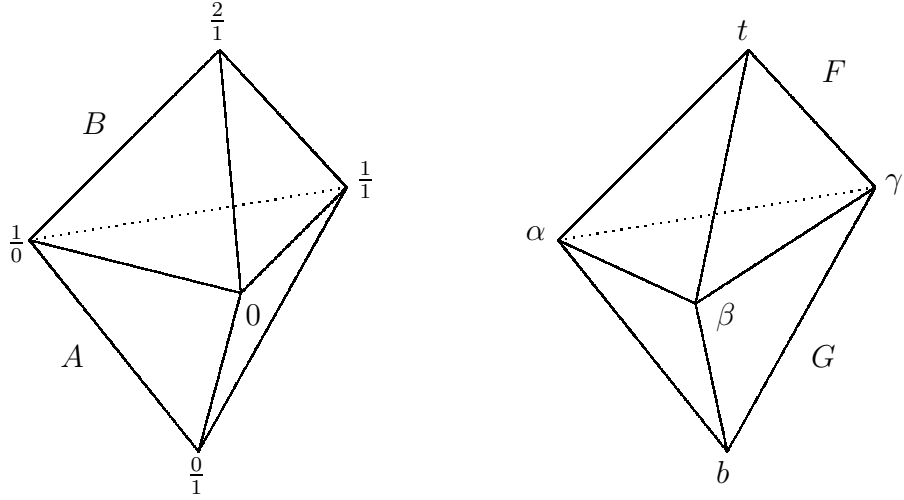
\begin{figure}[ht]
\begin{center}
\begin{picture}(100,42)

  \put(0,0){
    \qbezier(10,20)(20,17.5)(30,15)
    \qbezier(30,15)(35,20)(40,25)
    \pbez40(10,20)(25,22.5)(40,25)
    
    \qbezier(10,20)(10,20)(28,38)
    \qbezier(30,15)(30,15)(28,38)
    \qbezier(40,25)(40,25)(28,38)

    \qbezier(10,20)(10,20)(26,0)
    \qbezier(30,15)(30,15)(26,0)
    \qbezier(40,25)(40,25)(26,0)

\put(15,30){$B$}
\put(13,8){$A$}

\put(27,40){$\frac{2}{1}$}
\put(25,-3){$\frac{0}{1}$}
\put(41,25){$\frac{1}{1}$}
\put(8,19){$\frac{1}{0}$}
\put(30.5,12){$0$}

\put(50,0){
    \qbezier(10,20)(10,20)(23,14)
    \qbezier(23,14)(40,25)(40,25)
    \pbez40(10,20)(25,22.5)(40,25)
    
    \qbezier(10,20)(10,20)(28,38)
    \qbezier(23,14)(28,38)(28,38)
    \qbezier(40,25)(40,25)(28,38)

    \qbezier(10,20)(10,20)(26,0)
    \qbezier(23,14)(26,0)(26,0)
    \qbezier(40,25)(40,25)(26,0)

\put(35,35){$F$}
\put(34,8){$G$}

\put(27,39){$t$}
\put(25,-3){$b$}
\put(41,25){$\gamma$}
\put(7,20){$\alpha$}
\put(25,12){$\beta$}
    }
           } 

\end{picture}
\end{center}
\caption{The ideal triangulation of $M$}\label{f taucb}
\end{figure}
The tetrahedra labeled $A,B$ are the old ones of the complement of
the figure-eight knot. The gluing rules are the following:
$$
\begin{array}{lcr}
A^\frac{0}{1}\leftrightarrow B^\frac{2}{1}:\ 
(0,\frac{1}{0},\frac{1}{1})\leftrightarrow
(0,\frac{1}{0},\frac{1}{1})
&\qquad&
A^\frac{1}{0}\leftrightarrow B^0:\
(0,\frac{0}{1},\frac{1}{1})\leftrightarrow
(\frac{1}{0},\frac{1}{1},\frac{2}{1})\\

A^\frac{1}{1}\leftrightarrow B^\frac{1}{0}:\ 
(0,\frac{0}{1},\frac{1}{0})\leftrightarrow
(0,\frac{1}{1},\frac{2}{1})
&\qquad&
A^0\leftrightarrow F^\gamma:\
(\frac{0}{1},\frac{1}{1},\frac{1}{0})\leftrightarrow
(t,\alpha,\beta)\\

B^\frac{1}{1}\leftrightarrow G^\gamma:\ 
(0,\frac{1}{0},\frac{2}{1})\leftrightarrow
(b,\beta,\alpha)
&\qquad&
F^t\leftrightarrow G^b:\
(\alpha,\beta,\gamma)
\leftrightarrow
(\alpha,\beta,\gamma)\\

F^\alpha\leftrightarrow G^\beta:\
(\beta,\gamma,t)
\leftrightarrow
(\gamma,\alpha,b)
&\qquad&
F^\beta\leftrightarrow G^\alpha:\
(\alpha,t,\gamma)
\leftrightarrow
(\gamma,b,\beta)
\end{array}
$$

We call the moduli $z_A, z_B$ are referred to the
edge $\overline{0 \frac{1}{1}}$ and $z_F, z_G$ to $\overline{\alpha\beta}$.
The triangulation of the boundary torus is the one of 
Figure~\ref{f eqnsp}.

\setlength{\unitlength}{4pt}
\begin{figure}[ht]
\begin{center}
\begin{picture}(100,75)

  \put(0,25){\line(1,0){80}}
  \put(10,0){\line(1,0){60}}
  \multiput(20,25)(20,0){4}{\line(-2,-5){10}}
  \multiput(10,0)(20,0){4}{\line(-2,5){10}}
  \put(60,25){\line(2,5){10}}
  \put(80,25){\line(-2,5){10}}
  \put(20,25){\qbezier(0,0)(-15,40)(20,50)}
  \put(60,25){\qbezier(0,0)(15,40)(-20,50)}
  \put(20,25){\line(2,5){20}}
  \put(60,25){\line(-2,5){20}}
  \put(80,25){\line(2,5){10}}
  \put(70,50){\line(1,0){20}}
  \put(40,25){\line(2,5){10}}
  \put(40,25){\line(-2,5){10}}
  \put(30,50){\line(1,0){20}}

  \put(70,50){\line(2,5){10}}
  \put(90,50){\line(-2,5){10}}

  \multiput(0,0)(-10,25){2}{
    \put(10,0){\vector(1,0){10}}
    \put(77,17.5){\qbezier(0,0)(0,0)(0,-2)}
    \put(77,17.5){\qbezier(0,0)(0,0)(-1.5,-1.5)}

    \put(86,40){\qbezier(0,0)(0,0)(0,-2)}
    \put(86,40){\qbezier(0,0)(0,0)(-1.5,-1.5)}
    \put(85,37.5){\makebox(0,0){$*$}}
    }
  
  \multiput(0,0)(80,50){2}{
    \put(4,15){\line(0,-1){2.5}}
    \put(4,15){\qbezier(0,0)(0,0)(1.7,-1.8)}
    \put(4.6,11){$\bullet$}
    }
    
    \multiput(38,-0.75)(3,0){2}{$>$}
    \multiput(57,-0.75)(2,0){3}{$>$}
    
    \put(15.5,47){\qbezier(0,0)(0,0)(-1,-1.6)}
    \put(15.5,47){\qbezier(0,0)(0,0)(1,-1.6)}
    \put(15.6,49){\qbezier(0,0)(0,0)(-1,-1.6)}
    \put(15.6,49){\qbezier(0,0)(0,0)(1,-1.6)}

    \put(64.4,49){\qbezier(0,0)(0,0)(-1.2,1.4)}
    \put(64.4,49){\qbezier(0,0)(0,0)(0.8,1.7)}

    \put(64.2,51){\qbezier(0,0)(0,0)(-1.3,1.2)}
    \put(64.2,51){\qbezier(0,0)(0,0)(0.7,1.6)}

    \put(63.9,53){\qbezier(0,0)(0,0)(-1.4,1.2)}
    \put(63.9,53){\qbezier(0,0)(0,0)(0.6,1.6)}

%%%%%%%%%%%   LABELS %%%%%%%%%%%%%%%%%%%%%%%%
{\footnotesize

  \put(10,15){\makebox(0,0){$B_0$}}
  \put(10,23){\makebox(0,0){$\frac{1}{0}$}}
  \put(4,18){\makebox(0,0){$\frac{1}{1}$}}
  \put(14,14){\makebox(0,0){$\frac{2}{1}$}}

  \put(17,23){\makebox(0,0){$z_B$}}
  \put(4,23){\makebox(0,0){$\frac{1}{1-z_B}$}}
{\tiny  \put(10,8){\makebox(0,0){$1\!-\!\frac{1}{z_B}$}}}

  \put(20,10){\makebox(0,0){$A_0$}}
  \put(20,2){\makebox(0,0){$\frac{1}{1}$}}
  \put(16,11){\makebox(0,0){$\frac{0}{1}$}}
  \put(24,11){\makebox(0,0){$\frac{1}{0}$}}

  \put(20.25,20){\makebox(0,0){$z_A$}}

  \put(30,15){\makebox(0,0){$B_\frac{1}{0}$}}
  \put(30,23){\makebox(0,0){$\frac{1}{1}$}}
  \put(26,14){\makebox(0,0){$0$}}
  \put(34,14){\makebox(0,0){$\frac{2}{1}$}}

  \put(23,23){\makebox(0,0){$z_B$}}

  \put(30,35){\makebox(0,0){$G_\beta$}}
  \put(30,27){\makebox(0,0){$\gamma$}}
  \put(26,36){\makebox(0,0){$b$}}
  \put(34,36){\makebox(0,0){$\alpha$}}

  \put(23,27){\makebox(0,0){$z_G$}}

  \put(24,55){\makebox(0,0){$F_\beta$}}
  \put(18,53){\makebox(0,0){$\gamma$}}
  \put(24,38){\makebox(0,0){$t$}}
  \put(33,61){\makebox(0,0){$\alpha$}}

  \put(20,32){\makebox(0,0){$z_F$}}
  \put(33,68){\makebox(0,0){$1\!-\!\frac{1}{z_F}$}}

  \put(56,55){\makebox(0,0){$F_\alpha$}}
  \put(62,53){\makebox(0,0){$\gamma$}}
  \put(56,38){\makebox(0,0){$t$}}
  \put(47.5,61){\makebox(0,0){$\beta$}}

  \put(60,32){\makebox(0,0){$z_F$}}
  \put(47,68){\makebox(0,0){$\frac{1}{1-z_F}$}}

  \put(40,10){\makebox(0,0){$A_\frac{1}{0}$}}
  \put(40,2){\makebox(0,0){$0$}}
  \put(36,11){\makebox(0,0){$\frac{0}{1}$}}
  \put(44,11){\makebox(0,0){$\frac{1}{1}$}}

  \put(47,1.5){\makebox(0,0){$z_A$}}

  \put(50,15){\makebox(0,0){$B_\frac{2}{1}$}}
  \put(50,23){\makebox(0,0){$\frac{1}{1}$}}
  \put(46,14){\makebox(0,0){$\frac{1}{0}$}}
  \put(54,14){\makebox(0,0){$0$}}

  \put(57,23){\makebox(0,0){$z_B$}}

  \put(60,10){\makebox(0,0){$A_\frac{1}{1}$}}
  \put(60,2){\makebox(0,0){$0$}}
  \put(56,11){\makebox(0,0){$\frac{1}{0}$}}
  \put(64,11){\makebox(0,0){$\frac{0}{1}$}}

  \put(60.25,20){\makebox(0,0){$z_A$}}

  \put(70,15){\makebox(0,0){$B_\frac{1}{1}$}}
  \put(70,23){\makebox(0,0){$\frac{1}{0}$}}
  \put(66,14){\makebox(0,0){$\frac{2}{1}$}}
  \put(74,14){\makebox(0,0){$0$}}

  \put(63,23){\makebox(0,0){$z_B$}}

  \put(40,40){\makebox(0,0){$F_\gamma$}}
  \put(40,48){\makebox(0,0){$t$}}
  \put(36,39){\makebox(0,0){$\beta$}}
  \put(44,39){\makebox(0,0){$\alpha$}}

  \put(40.25,30){\makebox(0,0){$z_F$}}

  \put(50,35){\makebox(0,0){$G_\alpha$}}
  \put(50,27){\makebox(0,0){$\gamma$}}
  \put(46,36){\makebox(0,0){$\beta$}}
  \put(54,36){\makebox(0,0){$b$}}

  \put(57,27){\makebox(0,0){$z_G$}}

  \put(70,35){\makebox(0,0){$A_\frac{0}{1}$}}
  \put(70,27){\makebox(0,0){$\frac{1}{1}$}}
  \put(66,36){\makebox(0,0){$\frac{1}{0}$}}
  \put(74,36){\makebox(0,0){$0$}}

  \put(77,27){\makebox(0,0){$z_A$}}

  \put(80,40){\makebox(0,0){$F_t$}}
  \put(80,48){\makebox(0,0){$\beta$}}
  \put(76,39){\makebox(0,0){$\gamma$}}
  \put(83,39){\makebox(0,0){$\alpha$}}

  \put(87,48){\makebox(0,0){$z_F$}}

  \put(80,60){\makebox(0,0){$G_b$}}
  \put(80,52){\makebox(0,0){$\alpha$}}
  \put(76,61){\makebox(0,0){$\beta$}}
  \put(83,61){\makebox(0,0){$\gamma$}}

  \put(73,52){\makebox(0,0){$z_G$}}

  \put(40,60){\makebox(0,0){$G_\gamma$}}
  \put(40,52){\makebox(0,0){$b$}}
  \put(36,61){\makebox(0,0){$\beta$}}
  \put(43,61){\makebox(0,0){$\alpha$}}

  \put(40,70){\makebox(0,0){$z_G$}}

}

\end{picture} 
\end{center} 
\caption{Triangulation with moduli of the boundary torus}
\label{f eqnsp}
\end{figure}

We see that the system of compatibility and completeness
equations is equivalent to the following one:
$$\begin{array}{lcl}
\left\{\begin{array}{l}
\displaystyle{\frac{1}{1-z_A}\cdot\frac{1}{z_B}\cdot\frac{z_F}{z_G}=1}\\
z_Gz_F=1
\end{array}\right.&
\qquad&
\left\{\begin{array}{l}
\displaystyle{\frac{(1-Z_A)^2}{z_A}\cdot\frac{z_B^2}{1-z_B}=1}\\
z_B(1-z_A)=1
\end{array}\right.
\end{array}$$

From that easily we get $z_G=z_F$ and $z_F^2=1$. Since we are looking
for non degenerate solutions, we have $z_F=z_G=-1$. Using that we get
$z_A=z_B$ and $$z_A^2-z_A+1=0$$
and then $z_A=z_B=\displaystyle{\frac{1\pm i\sqrt 3}{2}}$. That is,
the ideal tetrahedra $F$ and $G$ are flat but not degenerate, while
$A$ and $B$ are regular, exactly as in the 
complement of the figure-eight knot. We notice that 
the space obtained by gluing
together the geometric versions of the tetrahedra $A,B,F,G$ is not a
manifold.

\end{document}